\documentclass[12pt, oneside]{amsart}
\usepackage{epsfig, lscape, amssymb, graphicx}

\usepackage{charter}
\usepackage{mathrsfs}

\usepackage[left=1.5in, right=1in, top=1.2in, bottom=1.2in]{geometry}

\usepackage{graphicx}
\usepackage{epstopdf}
\theoremstyle{plain}
\newtheorem{lemma}{Lemma}[section]
\newtheorem{theorem}[lemma]{Theorem}
\newtheorem{claim}[lemma]{Claim}
\newtheorem{conjecture}[lemma]{Conjecture}
\newtheorem{corollary}[lemma]{Corollary}

\newtheorem{proposition}[lemma]{Proposition}
\newtheorem{question}[lemma]{Question}

\newtheorem*{subsurface theorem}{Main Theorem}
\newtheorem*{girth conjecture}{Girth Conjecture}

\newcommand{\girth}{\operatorname{girth}}
\newcommand{\genus}{\operatorname{genus}}

\newcommand{\rel}{\operatorname{rel}}

\newcommand{\m}{\mathcal{M}}

\newcommand{\nn}{\mathcal{N}}
\newcommand{\bd}{\partial}
\newcommand{\hh}{\widetilde{H}}

\newcommand{\dd}{\mathbf{D}}
\newcommand{\db}{\mathbf{D}_B}
\newcommand{\dw}{\mathbf{D}_W}
\newcommand{\dg}{\mathbf{D}_G}

\newcommand{\dpm}{\mathbf{D}^\pm}

\newcommand{\ec}[1] { | \chi (#1) | }
\newcommand{\di}[2]{ [#1 | #2] }
\newcommand{\dint}[1]{ \mathbf{\underline{#1}} }

\newcommand{\hi}{ \wedge }
\newcommand{\comment}[1]{}
\newcommand{\R} { {\mathbb R} }

\newcommand{\Z} { {\mathbb Z} }

\newcommand{\scc} { \mathcal{S} }
\newcommand{\rr} { \mathcal{R} }

\newcommand{\locallyminimize}[1]{\texttt{LOCALLY.MINIMIZE(#1)}}
\newcommand{\segregation}[1]{\texttt{SEGREGATION(#1)}}
\newcommand{\constructcurve}[1]{\texttt{CONSTRUCT.CURVE(#1)}}
\newcommand{\constructdisc}[1]{\texttt{CONSTRUCT.DISC(#1)}}
\newcommand{\GIRTH}[1]{\texttt{GIRTH(#1)}}

\newcommand{\ps}{\mathscr{P\!S}}

\newcommand{\skipline}{\vspace{12pt}}

\begin{document}
\pagenumbering{roman}
\pagestyle{plain}


\title{The Girth of a Heegaard Splitting}
\author{Christopher John Jerdonek}
\date{\today}


\begin{center}

\begin{Large}The Girth of a Heegaard Splitting\end{Large}\\
\skipline
By\\
\skipline
CHRISTOPHER JOHN JERDONEK,\\
A.B. (Harvard University) 1998,\\
\skipline
a DISSERTATION\\
\skipline
submitted in partial satisfaction of the requirements for the degree of\\
\skipline
DOCTOR OF PHILOSOPHY\\
\skipline
in\\
\skipline
MATHEMATICS\\
\skipline
in the\\
\skipline
OFFICE OF GRADUATE STUDIES\\
\skipline
of the\\
\skipline
UNIVERSITY OF CALIFORNIA,\\
\skipline
DAVIS.\\
\vspace{.5in}
Approved:\\
\vspace{.25in}
William P. Thurston\\
\vspace{-.1in}\rule{2.5in}{1pt}\\
\medskip
Michael Kapovich\\
\vspace{-.1in}\rule{2.5in}{1pt}\\
\medskip
Greg Kuperberg\\
\vspace{-.1in}\rule{2.5in}{1pt}\\
\medskip
Joel Hass\\
\vspace{-.1in}\rule{2.5in}{1pt}\\
\medskip
Abigail Thompson\\
\vspace{-.1in}\rule{2.5in}{1pt}\\
\vspace{.2in}
Committee in Charge\\
\skipline
2005\\

\end{center}


\renewcommand{\baselinestretch}{1.6}\small\normalsize


\newpage
\setcounter{tocdepth}{1}
\tableofcontents


\newpage
\section*{Abstract}

We construct simple curves from immersed 
curves in the setting of handlebodies and Heegaard splittings.
We define a measure of complexity we call \emph{girth} for 
closed curves in a handlebody.
We extend this complexity to Heegaard splittings and pose 
a conjecture about all Heegaard splittings.
We prove a test case of this conjecture.
Let $S$ be a compact surface embedded in the boundary of 
a handlebody $H$.  Then the minimum girth over all curves in $S$
can be achieved by a simple closed curve.
We also present algorithms to compute the girth of curves
and surfaces.


\newpage
\section*{Acknowledgments}
Thank you to the UC Davis Mathematics Department for generously
supporting me throughout my graduate years.  I am very thankful
for their support.

I extend many thanks to everyone in the UC Davis Mathematics 
Department.
In particular, I thank the system administrators Zach Johnson, 
Marianne Waage, and Bill Broadley for consistently maintaining 
high quality computer resources.  I thank Celia Davis and Nancy 
Davis for their patience and administrative help.

I thank Professors Jack Milton and Tom Sallee for their teaching
advice.

I thank Professors Joel Hass, Abby Thompson, and Greg Kuperberg for 
their helpful comments on my mathematical work.  
I thank Ian Agol for interesting discussions in my earlier years.
I especially thank Professor Misha Kapovich for his advice 
and availability in my final year.

Most of all I would like to thank my advisor Bill Thurston 
for being a tremendous and lasting inspiration to me, both within 
mathematics and without.


\newpage
\pagestyle{headings}
\pagenumbering{arabic}
\markright{}

\section{Introduction}
\label{introduction}
Let $H\cup H^\prime$ be a Heegaard splitting of a closed orientable 3-manifold.    Let $\mathcal{C}$ denote the collection of closed curves in $\bd H$ that bound a disc in $H^{\prime}$.  Does the set of algebraically simplest elements of $\mathcal{C}$ overlap the set of geometrically simplest elements?  Below we make this question precise and begin an answer to it.

Let $\gamma$ be a closed curve in a handlebody $H$.  Define the girth of $\gamma$ to be the smallest intersection number between a curve homotopic to $\gamma$ and an essential disc in $H$.  Girth is an algebraic measure of complexity.  We say $\gamma$ is \emph{primitive} if $\gamma$ does not represent a proper power of some element of $\pi_{1}(H)$.  In particular, primitive curves are nontrivial in $\pi_{1}(H)$.  We propose the following answer to the above question.

\begin{girth conjecture}
Let $H\cup H^\prime$ be a Heegaard splitting of a closed orientable 
3-manifold.  Let $\gamma$ be a closed curve in $\bd H$ that is
primitive in $H$ and contracts in $H^\prime$.  Then there is an essential simple closed curve $\gamma^{\prime}$ in $\bd H$ with girth no greater than that of $\gamma$ which contracts in $H^\prime$.
\end{girth conjecture}

For example, let $M$ be a closed hyperbolic manifold.  Let $H\cup H^{\prime}$ be a Heegaard splitting of $M$.  Let $\gamma$ be an immersed curve in $\bd H$ that bounds a disc in $H^{\prime}$.  Does $\bd H$ contain an embedded curve $\gamma^{\prime}$ that bounds a disc in $H^{\prime}$ and has girth no greater than that of $\gamma$?  The Girth Conjecture would say yes.

Commonly in 3-manifold topology we wish to construct embedded objects from immersed objects, and Heegaard splittings are a fundamental setting.  The Girth Conjecture is like a Poincare Conjecture for all 3-manifolds.  By Proposition~\ref{finite fundamental group}, the Girth Conjecture for 3-manifolds with finite fundamental group is equivalent to the spherical space form conjecture.

\begin{figure}
\epsfig{file=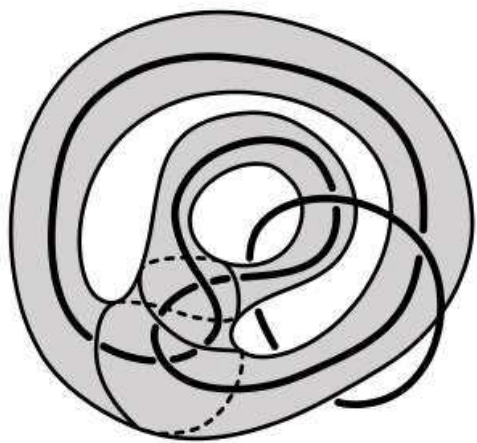, height=2.5in}
\caption{
All Dehn fillings of this standard Heegaard splitting $H\cup H^\prime$
of the Whitehead link $L$ satisfy 
Conjecture~\ref{ag girth conjecture}.
All hyperbolic manifolds that result have algebraic and
geometric girth four~\cite{Je2}. 
The thin curve bounds a disc in the gray handlebody $H$ that is disjoint from $L$.  This disc visibly intersects a disc in the complement $H^\prime$ 
in four points.}
\label{whitehead}
\end{figure}

Splittings of reducible 3-manifolds satisfy the Girth Conjecture by Lemma~\ref{reducible} and work of~\cite{Ha}.  In addition to spherical 3-manifolds, solvmanifolds satisfy the Girth Conjecture~\cite{Je1}, as well as the simplest genus two hyperbolic manifolds~\cite{Je2}.  See Figure~\ref{whitehead}.  The classification of Heegaard splittings of Seifert-fibered spaces~\cite{MS} provides an answer for genus two splittings of Seifert-fibered spaces.  We discuss connections to the geometry of 3-manifolds in more detail in Section~\ref{heegaardgirth}.

In Section~\ref{heegaardgirth}, we define the algebraic and geometric girths of a Heegaard splitting.  We show in Proposition~\ref{algebraic girth one} that the Girth Conjecture holds for all Heegaard splittings of algebraic girth one.  The proof uses the Freiheitssatz and Kneser's Conjecture for 3-manifolds.

Our Main Theorem shows that the Girth Conjecture is true in a certain general setting.  It is a purely topological result.

\begin{subsurface theorem}[\textbf{Theorem~\ref{subsurface theorem}}]
Let $H$ be a handlebody.  Let $S$ be a compact surface in 
$\partial H$.  The minimum girth in $H$ over all essential closed
curves in $S$ can be achieved by a simple closed curve in $S$. 
\end{subsurface theorem}

The Main Theorem implies the Girth Conjecture for curves $\gamma$ that are contained in an embedded ball in $H^\prime$.  If $B$ is a closed ball in $H^{\prime}$ that contains $\gamma$, simply take $S$ to be $B\cap \bd H$.  All curves in $B\cap \bd H$ contract in $H^{\prime}$.  Theorem~\ref{subsurface theorem} is trivial for surfaces $S$ that are compressible in $H$.  More interesting and typical surfaces $S$ can be obtained by taking the complement of a disjoint collection of simple closed curves in $\bd H$.

Our Main Theorem has the following corollaries.
Let $H$ be a handlebody.
A \emph{disc system} for $H$ is a disjoint collection of essential 
discs in $H$ that divide $H$ into a ball.
A Heegaard splitting $H\cup H^\prime$ is \emph{weakly reducible}
if there are essential discs in $H$ and $H^\prime$ whose
boundaries do not intersect in $\bd H$.  By ~\cite{CG}, weak
reducibility implies that either the 3-manifold is Haken,
or the splitting is reducible.

\begin{corollary}
\label{cor 1}
Let $H\cup H^\prime$ be a Heegaard splitting of a closed orientable 
3-manifold $M$.
Let $\delta$ be a collection of simple closed curves in $\partial H$ 
that bounds a disc system in $H^\prime$.
If there exists a curve in $\partial H-\delta$ that is homotopic in 
$H$ to a curve disjoint from some essential disc in $H$, 
then either $M$ is Haken, or the splitting is reducible.
\end{corollary}

\begin{proof}
Let $\gamma$ be a curve in $\partial H-\delta$ that is
homotopic in $H$ to a curve disjoint from an essential disc in $H$.
Then $\gamma$ has girth zero in $H$.  By the Main Theorem,
$\partial H-\delta$ contains a simple closed curve of girth zero
in $H$.

So we can assume $\gamma$ is simple.  By 
Lemma~\ref{one embedded curve}, the curve $\gamma$ is disjoint from 
an essential disc $D$ in $H$.
Since $\gamma$ is disjoint from a disc system in $H^\prime$,
$\gamma$ bounds an embedded disc $D^\prime$ in $H^\prime$.

Since $D$ and $D^\prime$ are essential discs in $H$ and $H^\prime$
with disjoint boundaries, the splitting $H\cup H^\prime$ 
is weakly reducible.  By~\cite{CG}, either 
the manifold is Haken, or the Heegaard splitting is reducible.
\end{proof}

The distance of a Heegaard splitting is one measure of 
complexity for Heegaard splittings.  Distance is explored,
for instance, in~\cite{He2}.
A distance one splitting is the same as a weakly reducible splitting.
Distance two also has special significance.
A Heegaard splitting $H\cup H^\prime$ has distance two if 
there exists a curve $\gamma$ in $\bd H$ that is disjoint from an 
essential disc in both $H$ and $H^\prime$.
This is also called the 
\emph{disjoint curve property}~\cite{Tho}.
 
\begin{corollary}
Let $H\cup H^\prime$ be a Heegaard splitting of a closed orientable 
3-manifold.
Let $D^\prime$ be an essential disc in $H^\prime$.
If there is a curve in $\partial H-\bd D^\prime$ that is homotopic 
in $H$ to a curve disjoint from some essential disc in $H$, then 
$H\cup H^\prime$ satisfies the disjoint curve property.
\end{corollary}

\begin{proof}
Let $\gamma$ be a curve in $\partial H-\bd D^\prime$ that is
homotopic in $H$ to a curve disjoint from some essential disc.
Then $\gamma$ has girth zero in $H$.  By the Main Theorem,
$\bd H-\bd D^\prime$ contains a simple closed curve of girth zero.

So assume that $\gamma$ is simple.
By Lemma~\ref{one embedded curve}, the curve $\gamma$ is disjoint 
from some essential disc $D$ in $H$.
Since the curve $\gamma$ is disjoint from both $D$ and $D^\prime$,
the splitting $H\cup H^\prime$ has distance two.
\end{proof}

We suggest two approaches to the Girth Conjecture for further
study.  Consider a Heegaard splitting $H\cup H^\prime$.
Let $\gamma$ be an immersed curve in $\bd H$ 
that bounds a disc $D^\prime$ in $H^\prime$.  Can we make 
$\gamma$ simple without increasing its girth in $H$?
The tower construction is one method of producing an 
embedded disc from an immersed curve.  Papakyriakapolous
used this method to prove the Loop Theorem.
See~\cite[p.45]{Hat} or~\cite[p.7]{Ja} for a description
of this construction.

If a neighborhood of the disc $D^\prime$ is a ball, we can 
apply the Main Theorem to $\gamma$.  If not, there 
are nontrivial coverings of this neighborhood.  The tower construction
proceeds by factoring the disc map through 
a sequence of 2-sheeted coverings.
This inspires the following question.

\begin{question}
Let $H\cup H^\prime$ be a Heegaard splitting of a closed $3$-manifold.
Let $S$ be a surface immersed in $\partial H^\prime$ so that $\pi_1(S)$
maps trivially into $H^\prime$.  Consider the minimum girth in $H$ 
over all essential closed curves in $S$.  Can the minimum be 
achieved by a simple closed curve in $S$?
\end{question} 

For embedded $S$ the answer is yes by the Main Theorem.
For $S$ an annulus the answer is yes by Lemma~\ref{power}.  
The answer is no if $S$ does not necessarily contract in $H^{\prime}$.  
To answer this question it suffices to consider the case of $S$
a thrice-punctured sphere.

Here is a second approach.
Let $S$ be a compact surface with boundary $\bd S$.
Let $\nn(\bd S)$ denote the subgroup of $\pi_1(S)$ normally 
generated by the boundary curves in $\bd S$.  This 
parallels the definition of the normal subgroup of relators 
in a Heegaard diagram.
The following statement resembles both the Girth Conjecture 
and the Main Theorem.

\begin{conjecture}
Let $H$ be a handlebody.  Let $S$ be a compact surface in $\bd H$.
Then the minimum girth in $H$ over all closed curves 
in $\nn(\bd S)$ can be achieved by a simple closed curve in $\nn(\bd S)$.
\end{conjecture}

In Section~\ref{trees} we introduce some basic results and constructions related to trees and handlebodies, and in Section~\ref{free groups} we do the same for free groups.  In Section~\ref{handlebodygirth} we present some elementary properties of girth.  We focus on simple closed curves.  There we also introduce a special conjecture, Conjecture~\ref{special conjecture}. In Section~\ref{heegaardgirth} we discuss girth for Heegaard splittings.  We define the algebraic and geometric girths of 
a Heegaard splitting and present a statement equivalent to the 
Girth Conjecture.  We also describe some results relating 
the girth of Heegaard splittings to 3-manifold geometry.

In Section~\ref{subordination} we prove Lemma~\ref{subordination-lemma}, an important lemma that we use several times throughout the paper.
In Section~\ref{curvealgorithm} we present an algorithm 
to compute the girth of a curve.  The algorithm uses Lemma~\ref{subordination-lemma}.

In Sections~\ref{maintheorem} and~\ref{mainlemma} we prove 
the Main Theorem.  We present several lemmas and the 
main proof in Section~\ref{maintheorem}.  We reserve the 
key Lemma~\ref{main lemma} for Section~\ref{mainlemma}.
The proof of Lemma~\ref{main lemma} contains the central construction 
used in the proof of the Subsurface Theorem.
Definitions needed to understand the proofs are spread throughout 
Sections~\ref{subordination} through~\ref{mainlemma}.

We conclude with an algorithm in Section~\ref{surfacealgorithm} to compute the girth for surfaces.  This short section pulls together results from several sections of the paper.

\section{Handlebodies, Trees, and Girth}
\label{trees}
In this section we introduce girth for closed curves in a handlebody.  We define girth and related concepts and prove some preliminary lemmas.

We use the language of ends.  See~\cite{E} or~\cite{St} for an algebraic introduction to the theory of ends.  See~\cite{GM} for a topological definition.  Let $X$ be a topological space with more than one end.  Define an \emph{interface} for $X$ to be a partition of the ends of $X$ into three subsets labeled black, white, and gray.  We require the black and white subsets of an interface to be nonempty.

A \emph{proper curve} in $X$ is a proper map of the real line to $X$.  Let $\gamma$ be a proper curve in $X$.  The proper curve $\gamma$ determines at most two ends of $X$.  These ends remain invariant under proper homotopies of $\gamma$ (see~\cite[p. 658]{GM}).  Let $\alpha$ be an interface for $X$.  We say $\gamma$ \emph{crosses} $\alpha$ if one end of $\gamma$ in $X$ is black and one is white, as determined by $\alpha$.

Let $H$ be a handlebody.  Let $D$ be an embedded essential disc in $H$.  The disc $D$ defines an interface for the universal cover of $H$, up to deck transformations, as follows.  Let $\hh$ denote the universal cover of $H$.  Consider a lift of $D$ to $\hh$.  This lift separates the ends of $\hh$ into two subsets we call black and white.  Define the gray subset to be empty.  This defines an interface $\alpha$ for $\hh$.  It is well-defined up to deck transformations.  We call $\alpha$ a \emph{lift} of $D$ to an interface for $\hh$.

Let $\gamma$ be a closed curve in $H$.  Let $D$ be an essential disc in $H$.  Define the homotopic intersection number of $\gamma$ and $D$ as follows.  Let $\alpha$ be a lift of $D$ to an interface for $\hh$.  The \emph{homotopic intersection number} of $\gamma$ and $D$ is the number of lifts of $\gamma$ that cross $\alpha$ in the universal cover of $H$.  We denote this quantity $\gamma\wedge D$.  By equivariance, this notion does not depend on the choice of lift $\alpha$.

Homotopic intersection number is invariant under homotopies of $\gamma$.  This is because homotopies of $\gamma$ lift to proper homotopies in the universal cover.  By Lemma~\ref{wedge} below, the homotopic intersection number of $\gamma$ and $D$ equals the minimum geometric intersection number of a curve homotopic to $\gamma$ and $D$.

The \emph{girth} of $\gamma$ in $H$ is the minimum of $\gamma\wedge D$ over all embedded essential discs $D$ in $H$.  We sometimes write this as $\girth(\gamma)$.  An essential disc $D$ is said to \emph{realize} the girth, or be a \emph{girth-realizing} disc, if $\gamma\wedge D=\girth(\gamma)$.  By Lemma~\ref{wedge} below, the girth of $\gamma$ equals the minimum geometric intersection number of an essential disc $D$ in $H$ and a curve homotopic to $\gamma$.

Let $H$ be a handlebody.  Let $\dd$ be a disjoint collection of essential discs in $H$.  Define the homotopic intersection number of $\gamma$ and $\dd$ to be the sum of $\gamma\hi D$ over all discs $D$ in $\dd$.   We denote this $\gamma\hi\dd$.

\begin{lemma}
\label{wedge}
Let $H$ be a handlebody.  Let $\gamma$ be a closed curve in $H$.  Let $\dd$ be a disjoint collection of essential discs in $H$.  The homotopic intersection $\gamma\wedge\dd$ equals the minimum geometric intersection number of $\dd$ and a curve homotopic to $\gamma$.
\end{lemma}

We prove Lemma~\ref{wedge} at the end of this section.  An immediate corollary of Lemma~\ref{wedge} is that if a closed curve $\gamma$ minimizes geometric intersection number with a disc collection $\dd$, then it also minimizes geometric intersection number with any subcollection of $\dd$:

\begin{corollary}
\label{subcollection}
Let $H$ be a handlebody.  Let $\gamma$ be a closed curve in $H$.  Let $\dd$ be a disjoint collection of essential discs in $H$.  Let $\dd_{1}$ be a subcollection of $\dd$.  If $\gamma$ minimizes geometric intersection with $\dd$ over its homotopy class, then $\gamma$ also minimizes geometric intersection with $\dd_{1}$ over its homotopy class.
\end{corollary}

\begin{proof}
Assume $\gamma$ does not minimize intersection with $\dd_{1}$.  It suffices to show that $\gamma$ does not minimize intersection with $\dd$.

Let $\dd_{2}$ denote the collection $\dd-\dd_{1}$.  Let $\gamma\hi\dd_{1}$ equal $n_{1}$, and let $\gamma\hi\dd_{2}$ equal $n_{2}$.  Since $\gamma$ does not minimize intersection with $\dd_{1}$, we have by Lemma~\ref{wedge} that $|\gamma\cap\dd_{1}|>n_{1}$.  Since $\gamma\hi\dd_{2}=n_{2}$, we also have $|\gamma\cap\dd_{2}|\geq n_{2}$, again by Lemma~\ref{wedge}.  Thus $|\gamma\cap\dd|>n_{1}+n_{2}$.

By the definition of homotopic intersection number for a collection of discs, the quantity $\gamma\hi\dd$ equals $n_{1}+n_{2}$.  By Lemma~\ref{wedge}, the minimum geometric intersection number between $\gamma$ and $\dd$ equals $n_{1}+n_{2}$.  So $\gamma$ does not minimize intersection with $\dd$.
\end{proof}

The language of trees will be useful throughout this paper and for the proof of Lemma~\ref{wedge}.  We assume familiarity with notions like paths, geodesics, and backtrackings in trees.  For instance, a geodesic in a tree is a path in the tree with no backtrackings.  See~\cite[Ch. I.2]{Se} for definitions of these terms and an introduction to the theory of trees.

Let $H$ be a handlebody.  Let $\dd$ be a disjoint collection of essential discs in $H$.  We associate an infinite tree $\Gamma$ to $\dd$ in the following way.  Consider all lifts of $\dd$ to the universal cover $\hh$ of the handlebody $H$.  Take the connected components of $\hh-\dd$ to be the vertices of $\Gamma$.  Take the lifts of $\dd$ to $\hh$ to be the edges of $\Gamma$.  Each edge connects two vertices, according to which components of $\hh-\dd$ the lift borders.  We call $\Gamma$ the \emph{tree associated to $\dd$}.

If the discs in $\dd$ are oriented, the graph $\Gamma$ is a directed graph.  If the connected components of $H-\dd$ are simply-connected, then $\Gamma$ is locally finite.  Otherwise, $\Gamma$ is locally infinite.  Here is an important special case.

Let $H$ be a handlebody.  A \emph{disc system} for $H$ is a collection $\dd$ of essential discs in $H$ such that $H-\dd$ is a ball.  If the discs are oriented, we sometime call it an oriented disc system.  If $g$ is the genus of $H$, any disc system has $g$ elements.  An oriented disc system for $H$ naturally gives rise to a collection of free generators for $\pi_{1}(H)$.  In this case, the tree associated to $\dd$ is also the Cayley graph of $\pi_{1}(H)$ with respect to the corresponding generators.  We will sometimes apply Lemma~\ref{wedge} to disc systems.

Let $\gamma$ be a closed curve in $H$ that is transverse to $\dd$.  In this case, we can lift $\gamma$ to a path in $\Gamma$ as follows.  Lift $\gamma$ to $\hh$.  Then map $\gamma$ to $\Gamma$ according to which discs the lift of $\gamma$ crosses.  The resulting path in $\Gamma$ we call a \emph{lift} of $\gamma$ to $\Gamma$.  When $\gamma$ is disjoint from $\dd$, the curve $\gamma$ lifts to a vertex of $\Gamma$, which we view as a degenerate, or constant path.

By the following lemma, $\gamma$ lifts to a path in $\Gamma$ without backtrackings if and only if the curve $\gamma$ minimizes intersection with $\dd$ over its homotopy class.  A path without backtrackings is either an infinite geodesic path or a vertex of $\Gamma$.

\begin{lemma}
\label{geodesic}
Let $H$ be a handlebody.  Let $\dd$ be a disjoint collection of essential discs in $H$.  Let $\Gamma$ be the tree associated to $\dd$.  Let $\gamma$ be a closed curve in $H$ that is transverse to $\dd$.  Lift $\gamma$ to $\Gamma$.  The curve $\gamma$ in $H$ minimizes intersection with $\dd$ in its homotopy class if and only if $\gamma$ in $\Gamma$ has no backtracking.
\end{lemma}

\begin{figure}
\includegraphics[height=2in]{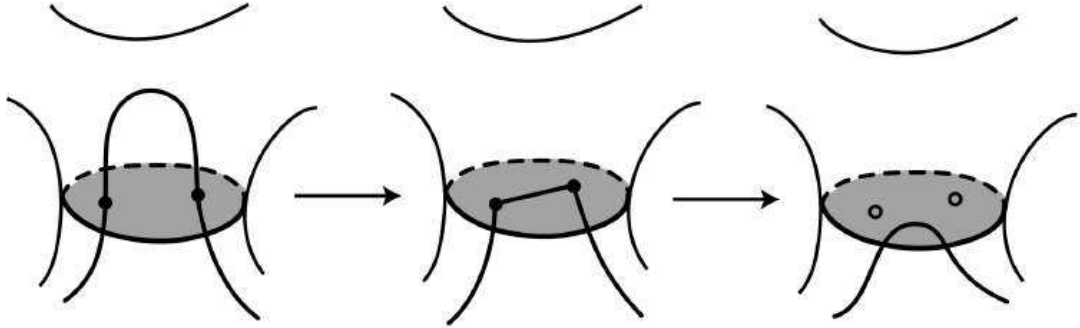}
\caption{The gray disc $D$ is part of a disjoint collection of discs $\dd$ in $H$.  Let $\Gamma$ be the tree associated to $\dd$.  If a lift of $\gamma$ to $\Gamma$ has a backtracking (shown above), then the intersection number of $\gamma$ with $\dd$ can be decreased by two in $H$.}
\label{homotopy 1}
\end{figure}

\begin{proof}
Assume first that $\gamma$ in $\Gamma$ has a backtracking.  We will show that $\gamma$ does not minimize intersection with $\dd$.  Choose a backtracking of $\gamma$ in $\Gamma$.  The backtracking corresponds to a component $\kappa$ of $\gamma-D$ in $H$ with the following two properties.  The closure $\overline{\kappa}$ of the arc $\kappa$ intersects $\dd$ only at its endpoints, and the arc $\kappa$ can be homotoped $\rel \bd \kappa$ into a disc in $\dd$.  Homotope $\kappa$ to be disjoint from $D$ in $H$, as shown in Figure~\ref{homotopy 1}.  This reduces the intersection of $\gamma$ and $\dd$ by two, as desired.

Assume now that $\gamma$ lifted to $\Gamma$ has no backtracking.  Then $\gamma$ lifts to either a vertex of $\Gamma$ or an infinite geodesic path in $\Gamma$.  In the former case, $\gamma$ intersects $\dd$ zero times, so it minimizes intersection with $\dd$.  Assume then that $\gamma$ lifts to an infinite geodesic path.  By the first paragraph, it suffices to show that this infinite geodesic is uniquely determined by the homotopy class of $\gamma$.

The infinite path $\gamma$ in $\Gamma$ defines two distinct ends of $\hh$.  Thus $\gamma$ defines two ends of $\Gamma$.   These two ends remain the same under homotopies of $\gamma$ in $H$.  In a tree, an infinite geodesic path is uniquely determined by its ends.  Since $\Gamma$ is a tree, the infinite geodesic determined by $\gamma$ is unique.
\end{proof}

We now prove Lemma~\ref{wedge}.

\begin{proof}[Proof of Lemma~\ref{wedge}]
Let $\gamma$ be any curve homotopic to the original $\gamma$.  Then $\gamma\wedge \dd$ is unchanged.  Consider any disc $D$ in $\dd$.  Lift $D$ to the universal cover $\hh$.  The disc $D$ defines an interface $\alpha$ for $\hh$.  If a lift of $\gamma$ crosses $\alpha$, then that lift must geometrically intersect the lift of $D$ at least once.  Thus $\gamma\wedge D \leq |\gamma\cap D|$.  Since this is true for all $D$, we have $\gamma\wedge\dd\leq|\gamma\cap \dd|$.

Now assume $\gamma\wedge \dd<|\gamma\cap\dd|$.  Then we have $\gamma\wedge D<|\gamma\cap D|$ for some $D$ in $\dd$.  It suffices to construct a homotopy of $\gamma$ that reduces the geometric intersection number of $\gamma$ and $\dd$.  Since $\gamma \wedge D<|\gamma\cap D|$, some lift of $\gamma$ must intersect $D$ at least twice.  In particular, $\gamma$ does not lift to a geodesic in $\Gamma$.  By Lemma~\ref{geodesic}, we can decrease the geometric intersection number of $\gamma$ and $\dd$.
\end{proof}

\section{Free Groups and Girth}
\label{free groups}
 We turn now to free groups.  Let $F$ be a free group.  Let $w$ be an element of $F$.  We define the girth of $w$ in $F$ as follows.  Let $g$ be the rank of $F$.  Form a handlebody $H$ of genus $g$.  Identify $F$ with the fundamental group $\pi_{1}(H)$.  Let $\gamma$ be a closed loop in $H$ that gives rise to $w$.  Define the girth of $w$ to be the girth of $\gamma$ in $H$.  We will sometimes write the girth of $w$ as $\girth(w)$.  The following lemma shows this notion is well-defined.  In other words, the girth of an element of a free group does not depend on the choice of identification between $F$ and $\pi_{1}(H)$.  A different identification differs by an automorphism of $\pi_{1}(H)$.

\begin{lemma}
Let $H$ be a handlebody.  Let $\gamma$ be a closed loop in $H$ that represents an element of $\pi_{1}(H)$.  Let $\Phi$ be an automorphism of $\pi_{1}(H)$.  Let $\gamma^{\prime}$ be a closed curve in $H$ that represents the element $\Phi(\gamma)$ in $\pi_{1}(H)$.  The girth of $\gamma$ in $H$ equals the girth of $\gamma^{\prime}$. 
\end{lemma}

\begin{proof}
The proof follows from the fact that automorphisms of a free group can be realized geometrically as homeomorphisms of a handlebody.

Let $F$ be a free group on $g$ letters.  The elementary Nielsen transformations of $F$ generate the group of automorphisms of $F$ (see for instance~\cite[Ch. I.4]{LSch}).  If $F$ is the fundamental group of a handlebody $H$ with generating system $\dd$, the elementary Nielsen transformations can be realized geometrically as a homeomorphism of $H$.  Thus all automorphisms of $F$ can be realized geometrically~\cite{Z}.

Let $n$ be the girth of $\gamma$.  Let $D$ be a girth-realizing disc for $\gamma$ in $H$.  By Lemma~\ref{wedge}, homotope $\gamma$ to intersect $D$ exactly $n$ times.  Realize the automorphism $\Phi$ of $\pi_{1}(H)$ as a homeomorphism of the handlebody $H$.

Apply the homeomorphism $\Phi$ to $H$, which contains $\gamma$ and $D$.  Then $\Phi(\gamma)$ is a curve homotopic to $\gamma^{\prime}$, and $\Phi(\gamma)$ intersects the disc $\Phi(D)$ exactly $n$ times.  Therefore the girth of $\gamma^{\prime}$ is less than or equal to $n$.  The reverse inequality follows by applying the same argument to the curve $\gamma^{\prime}$, in reverse.
\end{proof}

Let $F$ be a free group.  A \emph{generating system} for $F$ is a collection $\dd$ of free generators of $F$.  Let $g$ be the rank of $F$.  Since $F$ is Hopfian, any generating system for $F$ has $g$ elements.  We saw above that a handlebody $H$ with disc system $\dd$ gives rise to a generating system for $\pi_{1}(H)$.  

Let $\dd$ be a generating system for $F$.  Let $w$ be an element of $F$.  To parallel the definition above, we define $w\hi\dd$ as follows.  Let $H$ be a handlebody with disc system $\dd$.  Then $\pi_{1}(H)$ is naturally identified with $F$.  Choose a closed curve $\gamma$ that gives rise to $w$.  Define $w\hi\dd$ to equal $\gamma\hi\dd$.  In Section~\ref{computecurve}, we call $w\hi\dd$ the complexity of $\dd$ with respect to $w$.  The following lemma gives another meaning to $w\hi\dd$.

\begin{lemma}
\label{cyclic}
Let $F$ be a free group.  Let $\dd$ be a generating system for $F$.  Let $w$ be any element of $F$.  Then $w\hi\dd$ equals the length of $w$ written as a cyclically reduced word in the generators $\dd$ and their inverses.
\end{lemma}

Note that any element of the free group on $g$ letters has a unique cyclically reduced form, up to cyclic permutation (see~\cite[Sec. 1.4]{MKS}).  However, this is not needed for the proof.

\begin{proof}[Proof of Lemma~\ref{cyclic}]
Represent $F$ as the fundamental group of a handlebody $H$ with oriented disc system $\dd$.  The fundamental group of $H$ is in natural correspondence with $F$.  Consider any closed curve $\gamma$ in $H$.  The pattern of intersections with $\dd$ gives rise to a cyclic word $w$ in the generators $\dd$.  By Lemma~\ref{wedge}, the quantity $w\hi\dd$ equals the minimum intersection number of $\dd$ and a curve homotopic to $\gamma$.

It suffices to show that if the length of the cyclic word $w$ can be reduced, then $|\gamma\cap\dd|$ can be reduced.  Assume $w$ can be cyclically reduced in length.  By~\cite[Sec. 1.4]{MKS}, there is an adjacent pair of letters $ab$ in $w$ where $a=b^{-1}$.  Choose such a pair $ab$.  Consider the tree $\Gamma$ associated to $\dd$.  Lift $\gamma$ to $\Gamma$.  The pair $ab$ corresponds to a backtracking in the lift $\gamma$.  As in the proof of Lemma~\ref{wedge}, this backtracking gives rise to an elementary homotopy of $\gamma$ that decreases the intersection with $\dd$.  See Figure~\ref{homotopy 1} for the picture.
\end{proof}

Let $w$ be an element of the free group $F$.   By the following lemma, girth satisfies the group multiplicative property $\girth(w^{m})=m \girth(w)$.

\begin{lemma}
\label{power}
For $w$ in a free group $F$, $\girth(w^m)=m \girth(w)$.
\end{lemma}

This would follow from Lemma~\ref{cyclic} if girth-realizing discs were always non-separating.  We introduce the concept of a complete disc system for the proof.  A \emph{complete disc system} is a disjoint collection of essential discs in $H$ that cut $H$ into solid pairs of pants.  If $g$ is the genus of $H$, any complete disc system has $3g-3$ discs and cuts $H$ into $2g-2$ pairs of pants. 

\begin{proof}
Assume $w$ is nontrivial.  Let $g$ be the rank of $F$.  Form a handlebody $H$ of genus $g$.  Identify $F$ with the fundamental group $\pi_{1}(H)$.  Let $\gamma$ be a closed loop in $H$ that represents the element of $\pi_{1}(H)$ corresponding to $w$.  Then the girth of $w$ equals the girth of $\gamma$, and the girth of $w^{m}$ equals the girth of $\gamma^{m}$.

Let $n_{1}$ be the girth of $\gamma$.  Let $D_{1}$ be a girth-realizing disc for $\gamma$.  Freely homotope $\gamma$ to intersect $D_{1}$ exactly $n_{1}$ times.  Then $m$ copies of $\gamma$ intersect $D_{1}$ exactly $mn_{1}$ times.  Thus, $\girth(\gamma^m)\leq mn_{1}=m \girth(\gamma)$.

We now prove the reverse inequality.  Let $n_{2}$ be the girth of $\gamma^{m}$.  Let $D_{2}$ be a girth-realizing disc for $\gamma^m$.  Then $\gamma^{m}\hi D_{2}$ equals $n_{2}$.  Choose a complete disc system $\dd$ that contains $D_{2}$.  Freely homotope $\gamma^m$ to minimize intersection with $\dd$.  Since $\dd$ is a complete disc system and $\gamma^{m}$ is nontrivial, we have $|\gamma^{m}\cap\dd|>0$.  By Lemma~\ref{subcollection}, $\gamma^{m}$ also minimizes geometric intersection with $D_{2}$.  Furthermore, by Lemma~\ref{wedge}, the curve $\gamma^{m}$ intersects $D_{2}$ exactly $n_{2}$ times.

Let $\Gamma$ be the tree associated to $\dd$.  The graph $\Gamma$ is a trivalent graph.  By Lemma~\ref{geodesic}, the curve $\gamma^{m}$ lifts to a path $P$ in $\Gamma$ without backtrackings.  Since $\gamma^{m}$ is not disjoint from $\dd$, the path $P$ is an infinite geodesic path.  The action of $F$ on $\hh$ via deck transformations induces an action of $F$ on $\Gamma$.  The element $w^{m}$ leaves $P$ invariant.  Indeed, $P/\langle w^{m}\rangle$ is a path in $\Gamma/\langle w^{m}\rangle$ that crosses the edge corresponding to $D_{2}$ exactly $n_{2}$ times.

To complete the proof, it suffices to show that $w$ also leaves $P$ invariant.  For if this were the case, then $P/\langle w\rangle$ would be a circuit crossing the edge corresponding to $D_{2}$ exactly $n_{2}/m$ times.  The sequence of edges crossed by this circuit would define a loop in $H$ that is freely homotopic to $\gamma$ and that crosses $D_{2}$ exactly $n_{2}/m$ times.  Hence, $\girth(w)\leq n_{2}/m=\girth(\gamma^{m})/m$.

To see that $w$ leaves $P$ invariant, observe that $w^{m}$ leaves the path $w(P)$ invariant.  By a theorem of Tits~\cite[Prop. 24]{Se}, $w(P)$ must contain the infinite geodesic $P$.  Thus $w(P)$ and $P$ are actually the same.
\end{proof}

\section{Girth of Simple Closed Curves}
\label{handlebodygirth}
In this section we prove some basic properties about the girth of simple closed curves on the boundary of a handlebody.  For instance, we show that there are simple closed curves of arbitrarily high girth on the boundary of any handlebody.  We also introduce a conjecture equating the girth of a simple 
closed curve to another simple quantity.

Let $H$ be a handlebody.  Let $\gamma$ be a simple closed curve in $\bd H$.  
The main result of this section is the following Lemma~\ref{one embedded curve}.  When $\gamma$ is a simple closed curve, the minimum geometric intersection number of a curve homotopic to $\gamma$ and an essential disc in $H$ can be achieved by keeping $\gamma$ fixed.  In other words, no homotopy of $\gamma$ is necessary to realize the girth.

\begin{lemma}
\label{one embedded curve}
\label{simple closed curve}
Let $H$ be a handlebody.  Let $\gamma$ be a simple closed curve in $\bd H$.  Let $D$ be a girth-realizing disc for $\gamma$.  There exists an essential disc $D^{\prime}$ in $H$ such that $|\gamma\cap D^{\prime}|$ equals $\gamma\wedge D$ .  Equivalently, the girth of $\gamma$ equals the minimum geometric intersection number of $\gamma$ and an essential disc in $H$.
\end{lemma}

We give two proofs of Lemma~\ref{one embedded curve} in this section.  The first is somewhat shorter, but the construction is less direct.  It uses the tower construction in the proof of Dehn's Lemma.  The second proof is longer.  The construction is more explicit and easier to carry out.  We defer that proof to the end of this section.  There is even a third proof of Lemma~\ref{one embedded curve}.  The third proof can be obtained by following the proof of our main result, Theorem~\ref{subsurface theorem} in Section~\ref{maintheorem}.  That proof is of course much longer but entirely different in spirit.

\begin{proof}[First proof of Lemma~\ref{one embedded curve}.]
Let $n$ be the girth of $\gamma$.  The disc $D$ in $H$ is a girth-realizing disc for $\gamma$.  We will construct an essential disc $D^{\prime}$ in $H$ that intersects $\gamma$ exactly $n$ times.  Isotope $D$ in $H$ so that $D$ is transverse to $\gamma$.  Choose a single lift of $D$ to the universal cover $\pi:\hh\rightarrow H$.  The disc $D$ in the universal cover $\hh$ defines an interface $\beta$ for $\hh$.  By the definition of girth, the number of lifts of $\gamma$ that cross $\beta$ equals $n$. 

Let $B$ denote the boundary $\bd D$ of this lift to $\hh$.  The simple closed curve $B$ represents a nontrivial element of $H_{1}(\bd \hh,\Z)$, where $H_{i}$ denotes the usual homology.  We will construct a disjoint collection of simple closed curves in $\bd \hh$ that is homologous to $B$ and intersects $\pi^{-1}(\gamma)$ exactly $n$ times.  

Note that every curve in $\pi^{-1}(\gamma)$ can be homotoped to intersect $B$ either zero or one times.  The curves that cross $\beta$ can be made to intersect $B$ once, while the curves that do not cross $\beta$ can be made to intersect $B$ not at all.

View $B$ as a collection of curves.  Consider the lifts of $\gamma$ that intersect $B$ geometrically.  Since $|D\cap\gamma|$ is finite, there are finitely many such lifts.  For each such lift, modify $B$ in a small neighborhood so that the lift intersects the new $B$ either zero or one time.  This can be done while keeping the homology class of $B$ the same.  Lifts that cross $\beta$ will intersect $B$ once, while lifts that do not cross $\beta$ will be disjoint from $\beta$.  See Figure~\ref{curvemodify} for an illustration of how this can be done.

\begin{figure}
\epsfig{file=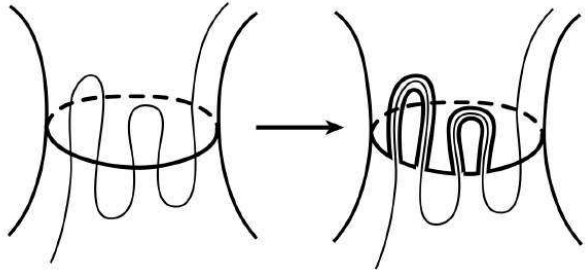, height=1.5in}
\caption{
This lift of $\gamma$ to the universal cover $\hh$ intersects the curve collection $B$ five times.  Since the lift crosses the interface $\beta$ once, the curve $B$ can be modifed to intersect the lift just once.  The change takes place in a neighborhood of the lift.  The new $B$ is homologous to the original.
}
\label{curvemodify}
\end{figure}

We now have a curve collection $B$ that intersects $\pi^{-1}(\gamma)$ exactly $n$ times.  Choose a curve component $B^{\prime}$ of $B$ that represents a nontrivial element of $H_{1}(\bd \hh,\Z)$.  This is possible since the original $B$ represented a nontrivial homology class.  This curve bounds a disc $D^{\prime}$ in $\hh$.

Project the disc $D^{\prime}$ to $H$.  The disc $D^{\prime}$ intersects $\gamma$ no more than $n$ times in $H$.  However, $D^{\prime}$ is not necessarily embedded.  To complete the proof, we apply the following Lemma~\ref{tower} to $D^{\prime}$.
\end{proof}

Lemma~\ref{tower} is where we use the tower construction.  Our other proof of Lemma~\ref{one embedded curve} is different.  In that proof, more care is taken in modifying $B$.  The modification is done in an equivariant way.  In all stages of that construction, the curves $B$ project down to a disjoint collection of simple closed curves in $H$.

\begin{lemma}
\label{tower}
Let $H$ be a handlebody.  Let $\gamma$ be a simple closed curve in $\bd H$.  Let $D$ be a disc in $H$ whose boundary is essential in $\bd H$ but is not necessarily embedded.  Then there is an essential embedded disc $D^{\prime}$ in $H$ with $|D^{\prime}\cap\gamma|\leq|D\cap\gamma|$.
\end{lemma}

For this proof, we assume familiarity with the tower construction.  In the course of the proof, observe that Lemma~\ref{tower} continues to hold true if the handlebody $H$ is replaced by any compact 3-manifold.  Moreover, if $\bd H$ has a metric, the result holds true if $|D\cap\gamma|$ is replaced by measuring the length of $\bd D$.  Note also that $\gamma$ need not be simple nor connected for the proof to hold.  For all of these variations, the same argument works.

\begin{proof}[Proof of Lemma~\ref{tower}.]
Assume that $\bd D$ is not a simple closed curve.  Homotope $\gamma$ in $\bd H$ so that $\gamma$ is transverse to $D$ and so that $\gamma\cap D$ lies away from the self-intersections of $\bd D$.  This can be done without increasing the intersection number $|\gamma\cap D|$.  Also make the self-intersections of $\bd D$ transverse, again without increasing $|\gamma\cap D|$.

Since we will be modifying $D$, view $D$ as a map $D^{2}\rightarrow D_{0}\subset H$.  Consider a regular neighborhood $V_{0}$ of $D_{0}$ in $H$.  In the tower construction, one considers a finite sequence of 2-sheeted covers, starting with a 2-sheeted cover of $V_{0}$.  At each stage of the construction, we lift the map of $D^{2}$.  At the top of the tower construction, we have a space $V_{n}$ and a map $D^{2}\rightarrow D_{n}\subset V_{n}$.  The space $V_{n}$ has no 2-sheeted cover, so $\bd V_{n}$ consists of 2-spheres.

Let $S^{2}$ denote the component of $\bd V_{n}$ that contains $\bd D_{n}$.  Consider the preimage of $\bd H$ in $S^{2}$ under the composition of the covering maps.  This is a compact subsurface of $S^{2}$ that contains $\bd D_{n}$.  Call this compact subsurface $S$.

At this point we need a simple closed curve $C$ in $S$ that projects to an essential closed curve in $\bd H$ and does not increase intersection with $\gamma$.  There are a couple ways to do this.  We choose to apply Lemma~\ref{top of tower}.  Rather than interrupt the flow of the proof, we prove this lemma afterwards.




Apply Lemma~\ref{top of tower} to $\bd D_{n}$ mapping to $\pi_{1}(\bd H)$.  We obtain a simple closed curve $C$ in $S$ that projects to an essential closed curve in $\bd H$.  Since $C$ is a subset of $\gamma$, it does not increase intersection with $\gamma$ downstairs.  Thus, where $\pi$ is the composition of the projection maps down to $H$, we have $|\pi(C)\cap \gamma|\leq|D\cap\gamma|$.

Since the simple closed curve $C$ lies inside the 2-sphere $S^{2}$, it bounds an embedded disc $B_{n}$ inside $V_{n}$.  For the remainder of the proof, we descend through our tower of 2-sheeted covers, modifying $B_{n}$ at each stage.

At each stage, map the embedded disc $B_{i}$ down to a disc $B_{i-1}$ via the corresponding covering map.  Since the cover is 2-sheeted, $B_{i-1}$ has at most double points.  These can be resolved in the standard way by cutting-and-pasting.  This cutting-and-pasting does not increase intersection with $\gamma$ in $H$ because $\bd B_{i-1}$ is affected only near the preimage of self-intersections of $\bd D$.  The intersection points $D\cap\gamma$ were taken to be away from the self-intersections of $\bd D$.  In the end, we have an embedded disc $B_{0}$ in $H$.  Take $D^{\prime}$ to be $B_{0}$.  $D^{\prime}$ is an embedded essential disc in $H$ satisfying $|D^{\prime}\cap\gamma|\leq|D\cap\gamma|$.
\end{proof}

We now prove the final lemma needed for the proof of Lemma~\ref{one embedded curve}.  Lemma~\ref{top of tower} is needed at the top of the tower construction in the proof of Lemma~\ref{tower}.

\begin{lemma}
\label{top of tower}
Let $S$ be a surface.  Let $G$ be a group.  Let $\Phi:\pi_{1}(S)\rightarrow G$ be a homomorphism from the fundamental group of $S$ to $G$.  Let $\gamma$ be a curve in $S$ with transverse self-intersections.  Suppose $\Phi$ maps $\gamma$ to a nontrivial conjugacy class in $G$.  Then $\gamma\subset S$ contains a simple closed curve that maps to a nontrivial conjugacy class in $G$.
\end{lemma}

\begin{figure}
\epsfig{file=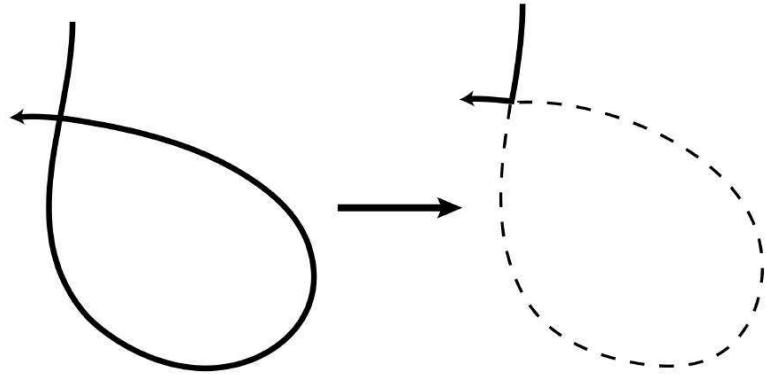, height=2in}
\caption{
We show part of a curve $\gamma$ in a surface $S$.  The curve $\gamma$ maps to a nontrivial conjugacy class in a group $G$.  If the loop on the left maps trivially, then removing it creates a curve with fewer self-intersections that continues to map nontrivially.
}
\label{selfintersection}
\end{figure}

\begin{proof}
The proof is by downward induction on the number of self-intersections of $\gamma$.  Choose a point in $\gamma$ and begin tracing out the curve $\gamma$ in one direction.  Stop when the curve first intersects itself.

After dropping the initial arc segment, this forms a simple closed curve $\gamma^{\prime}$ in $S$.  if $\gamma^{\prime}$ maps to a nontrivial conjugacy class in $G$, we are done.  If $\gamma^{\prime}$ maps trivially to $G$, then delete $\gamma^{\prime}$ from $\gamma$ to form a new curve $\gamma$.  See Figure~\ref{selfintersection} for an illustration.  The homomorphism $\Phi$ still maps the new $\gamma$ to a nontrivial conjugacy class.

Repeat the construction until the proof is complete.
\end{proof}

We now return to our main narrative.  Let $\pi:\hh\rightarrow H$ denote the universal covering space of $H$.  Define a graph $\Gamma$ in $\bd \hh$ as follows.  Take the components of $\bd\hh-\pi^{-1}(\gamma)$ to be the vertices of $\Gamma$.  Take the components of $\pi^{-1}(\gamma)$ to be the edges.  Each component of $\pi^{-1}(\gamma)$ is a curve that borders two components of $\bd\hh-\pi^{-1}(\gamma)$, which are possibly the same.  In this way, each edge defines two vertices.  We call $\Gamma$ the \emph{graph associated} to $\gamma$.  The graph $\Gamma$ is locally infinite if the genus of $H$ is two or more.

The girth of a graph is the length of a shortest cycle in the graph.  By the following lemma, the girth of a simple closed curve $\gamma$ equals the girth of the graph $\Gamma$ associated to $\gamma$.  This explains our choice of the word girth.  Lemma~\ref{graph girth} is a corollary of Lemmas~\ref{one embedded curve} and~\ref{tower}.

\begin{lemma}
\label{graph girth}
Let $H$ be a handlebody.  Let $\gamma$ be a simple closed curve in $\bd H$.  Let $\Gamma$ be the graph associated to $\gamma$.  If $\bd H-\gamma$ is incompressible in $H$, then the girth of $\gamma$ equals the girth of $\Gamma$.
\end{lemma}

By Lemma~\ref{one embedded curve}, the surface $\bd H-\gamma$ is compressible if and only if $\gamma$ has zero girth.  

\begin{proof}
Let $n$ be the girth of $\gamma$.  By Lemma~\ref{one embedded curve}, the smallest intersection number of $\gamma$ with an essential disc in $H$ equals $n$.  Let $D$ be an essential embedded disc that achieves this minimum.  Then $D$ intersects $\gamma$ in $n$ points.  Since $\bd H-\gamma$ is incompressible, $n$ is bigger than zero.

Lift $D$ to the universal cover $\hh$.  Since $n$ is bigger than zero, the boundary of $D$ defines a circuit $C$ in the graph $\Gamma$ associated to $\gamma$.  The circuit $C$ has $n$ edges.  Thus, the girth of $\Gamma$ is less than or equal to the girth of $\gamma$.

We now prove the reverse inequality.  Let $m$ be the girth of the graph $\Gamma$ in the graph theory sense.  Let $C$ be a circuit in $\Gamma$ that achieves the girth of $\Gamma$.  Then the length of $C$ is $m$.  Realize $C$ as the boundary of a disc $D$ in $\hh$.  Then $D$ intersects lifts of $\gamma$ exactly $m$ times. Map $D$ to $H$.  Then $D$ is a disc in $H$ that intersects $\gamma$ exactly $m$ times but is not necessarily embedded.

Apply Lemma~\ref{tower} to the non-embedded disc $D$ to create an embedded disc $D$ that intersects $\gamma$ no more than $m$ times.  
The girth of $\gamma$ equals the minimum geometric intersection number of an essential embedded disc in $H$ and a curve homotopic to $\gamma$.  Thus the girth of $\gamma$ is less than or equal to $m$, the girth of $\Gamma$.
\end{proof}

In the case that $\gamma$ is simple, we conjecture another interpretation of girth.  We conjecture that the girth of $\gamma$ equals the minimum geometric intersection number of $\gamma$ and a curve in $\bd H$ that is homotopic to $\gamma$ in $H$ but not in $\bd H$.  In other words, the girth of $\gamma$ could be computed by ranging over all curves homotopic to $\gamma$ rather than over all essential discs in $H$.

\begin{conjecture}
\label{special conjecture}
Let $H$ be a handlebody.  Let $\gamma$ be a simple closed curve 
in $\bd H$.  Let $\kappa$ be a closed curve in $\bd H$ 
that is homotopic to $\gamma$ in $H$ but not isotopic to 
$\gamma$ in $\bd H$.  Then the intersection number of $\kappa$ and 
$\gamma$ is at least the girth of $\gamma$.
\end{conjecture}

To obtain a curve $\kappa$ homotopic to $\gamma$ such that $|\kappa \cap \gamma|$ equals the girth of $\gamma$, simply Dehn-twist $\gamma$ around a girth-realizing disc for $\gamma$. 

By Lemma~\ref{one embedded curve}, a simple closed curve of girth zero is a curve disjoint from some essential disc in $H$.  By the following lemma, simple closed curves of girth one do not exist in the higher genus case.  The argument is the classical argument that stabilized Heegaard splittings of genus bigger than one are reducible.  We repeat it here for completeness.
  
\begin{lemma}
\label{geometric girth one}
Let $H$ be a handlebody of genus two or more.  Let $\gamma$ be a closed curve in $\partial H$.  If $\gamma$ is simple, then $\gamma$ does not have girth one.
\end{lemma}

\begin{proof}
Let $\gamma$ be a curve in $\bd H$ satisfying $\girth(\gamma)\leq 1$ in $H$.
By Lemma~\ref{simple closed curve}, there is an essential disc $D$ that
intersects $\gamma$ either once or not at all.  Suppose a disc $D$ in $H$ transversely intersects $\gamma$ exactly once. Consider the boundary of a regular neighborhood in $\bd H$ of $\gamma\cup\bd D$.
This curve bounds a disc in $H$ and intersects $\gamma$ zero times.
Since $H$ has genus bigger than one, the curve is essential in $\bd H$.
Thus the girth of $\gamma$ is zero.
\end{proof}

Lemma~\ref{geometric girth one} is actually true for all curves $\gamma$ in $H$, but the proof is more difficult.  We prove this as Lemma~\ref{group girth one}.  For now we prove the following special case.  Lemma~\ref{genus two} is a corollary of Lemma~\ref{geometric girth one}.

\begin{lemma}
\label{genus two}
Let $H$ be a handlebody.  Let $\gamma$ be a closed curve in $H$.  If the genus of $H$ equals two, then $\gamma$ does not have girth one.
\end{lemma}

\begin{proof}
Let $\gamma$ be a closed curve in $H$ that intersects an essential disc $D$ exactly once.  It suffices to show that $\gamma$ has girth zero.  Extend $D$ to a disc system $\dd$ for $H$.  Let $\{a, b\}$ be the corresponding generating system for $\pi_{1}(H)$, with $a$ corresponding to $D$.  The curve $\gamma$ is conjugate to a group element of the form $a b^{m}$, for some integer $m$.  Homotope $\gamma$ to a simple closed curve in $\bd H$ that intersects $D$ exactly once.  By Lemma~\ref{geometric girth one}, the girth of $\gamma$ must be zero.
\end{proof}

Whereas Lemma~\ref{geometric girth one} shows that simple closed curves of girth one do not exist, Lemma~\ref{large girth} shows that simple closed curves of large girth do exist.  In particular, every free group contains elements of arbitrarily high girth.  

\begin{lemma}
\label{large girth}
Let $H$ be a handlebody.  The surface $\bd H$ contains simple closed
curves of arbitrarily high girth in $H$.
\end{lemma}

\begin{proof}
Let $g$ be the genus of $H$.  Let $S$ be a closed disc with $g$ open discs removed.  Then $H$ is homeomorphic to the product of $S$ with a closed
interval.  Identify $H$ with $S\times I$.  The boundary of the cross-section $S\times\{1/2\}$ is a collection of $g+1$ closed curves in $\bd H$.  Denote this collection $\gamma$. 

The collection of curves $\gamma$ represents a rational element of $\ps(\bd H)$, the projective space of measured laminations on $\bd H$.  See~\cite{Thu2} an understanding of $\ps(\bd H)$.  Perturb $\gamma$ to an irrational element $\lambda$ in $\ps(\bd H)$.  The simple closed curves in $\bd H$ are dense in $\ps(\bd H)$.  By approximating $\lambda$ by simple closed curves, we can generate curves $\gamma^\prime$ of higher and higher girth.

Consider a standard collection of $3g-3$ non-parallel essential discs in $H$, so that each disc intersects $\gamma$ twice.  See Figure~\ref{highgirth}.  Call this collection $\dd$.  The collection $\dd$ divides $H$ into $2g-2$ filled-in pairs of pants.  Each pair of pants has three boundary arcs that connect the three disc boundaries in a cycle.

\begin{figure}
\epsfig{file=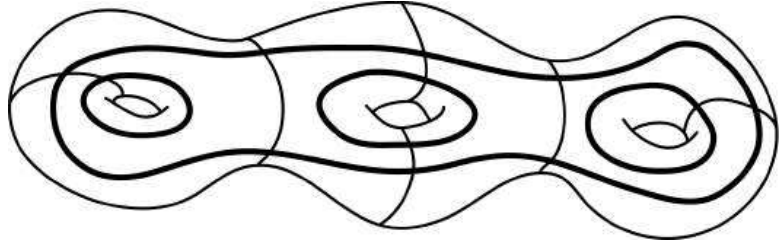, height=1.25in}
\caption{
The thin curves are the boundaries of $3g-3$ discs in a handlebody of genus $g$.  Each disc intersects the collection $\gamma$ of thick curves in exactly two points.
}
\label{highgirth}
\end{figure}

Consider the measured lamination $\lambda$ with respect to $\dd$.
If the perturbation is small enough, the lamination 
$\lambda$ will intersect each pair of pants in three bands of arcs 
in the same isotopy class as the original arcs.  See 
Figure~\ref{highgirthpant}.  Each band contains an infinite number 
of arcs, and each band has roughly the same measure.  

\begin{figure}
\epsfig{file=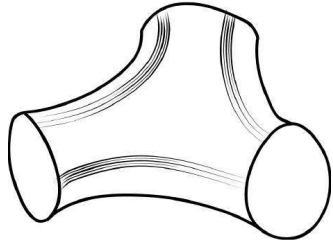, height=1.25in}
\caption{
Perturb the curve $\gamma$ of Figure~\ref{highgirth}  
to an irrational measured lamination $\lambda$.
The lamination $\lambda$ intersects each of the $2g-2$ pants 
in three bands of arcs, after isotopy.}
\label{highgirthpant}
\end{figure}

Let $\gamma^\prime$ be a simple curve that approximates $\lambda$.
The curve $\gamma^\prime$ must intersect each pair of pants in 
three bands of arcs in the same isotopy class as the bands for 
$\lambda$.  Each band for $\gamma^\prime$ has roughly the same 
number of arcs.  Let $n$ denote the smallest number of these 
arcs among all $6g-6$ bands of arcs.  As $\gamma^\prime$ more 
closely approximates $\lambda$, the number $n$ must grow larger.

Let $D$ be an essential disc in $H$ that minimizes intersection with $\gamma$.  By Lemma~\ref{one embedded curve}, the geometric intersection number of $D$ and $\gamma$ equals the girth of $\gamma$.  It is evident that $D$ must intersect $\gamma^\prime$ in at least $2n$ points.  See Figure~\ref{highgirthpant} for an illustration.  Thus the girth of $\gamma$ is at least $2n$.
\end{proof}

We now present our second proof of Lemma~\ref{one embedded curve}.

\begin{proof}[Second proof of Lemma~\ref{one embedded curve}.]
Let $n$ be the girth of $\gamma$.  Recall that $\gamma$ is a simple closed curve in $\bd H$.  Let $D$ be a disc in $H$ that realizes the girth of $\gamma$.  We will construct a disc $D^{\prime}$ in $H$ that intersects $\gamma$ exactly $n$ times.

Let $\pi:\hh\to H$ be the universal cover.  Let $\Gamma$ be the tree associated to $D$.  The action of $\pi_{1}(H)$ on $\hh$ induces an action of $\pi_{1}(H)$ on $\Gamma$.  Note that the action of $\pi_{1}(H)$ on $\Gamma$ is not free when $\genus(H)\geq 2$.  Isotope $\gamma$ in $\bd H$ to be transverse to $D$.  Lift $\gamma$ to a path $P$ in $\Gamma$.  The path $P$ is invariant under the cylic action of $\langle\gamma\rangle$ on $\Gamma$.

Descend to $\Gamma/\langle\gamma\rangle$.  Let $C$ denote the path $P/\langle\gamma\rangle$ in $\Gamma/\langle\gamma\rangle$.  The path $C$ is a map of a circuit to $\Gamma/\langle\gamma\rangle$.  The edges of $C$ are in correspondence with the points $\gamma\cap D$ in $H$, and the vertices of $C$ are in correspondence with the arcs $\gamma-D$ in $H$.

Consider a single lift of $D$ to $\hh$.  Let $B$ denote the boundary of this lift.  The curve $B$ is a simple closed curve that separates $\bd\hh$ into two noncompact pieces.  It represents a nontrivial element of $H_{1}(\bd \hh,\Z)$, where $H_{i}$ denotes the usual homology.  The curve $\pi(B)$ is the same as $\bd D$.  The points $\gamma\cap\pi(B)$ are in correspondence with $\gamma\cap D$.

Rename $\gamma$ as $\gamma^{\prime}$.  We reserve $\gamma$ for the original curve.  The proof proceeds by downward induction on the number of edges of $C$.    At each stage we modify the path $C$ in $\Gamma/\langle\gamma\rangle$, the curve $\gamma^{\prime}$ in $H$, and the curve $B$ in $\bd \hh$.  We stop when $C$ is a geodesic path in $\Gamma/\langle\gamma\rangle$.  Throughout the induction we preserve the correspondence between the edges of $C$, the points of $\gamma^{\prime}\cap D$, and the points of $B\cap \gamma$.  Note that after the first stage, the curve collection $B$ can have more than one connected component.  

\textbf{*Begin induction.}
Consider the path $C$.  If $C$ has no backtracking  in $\Gamma/\langle\gamma\rangle$, then skip to the end of the induction stage (\textbf{**}).  Otherwise, consider a backtracking of $C$.  The backtracking consists of a vertex and two edges of $C$.  The vertex of the backtracking corresponds to an arc $\kappa^{\prime}$ of $\gamma^{\prime}-D$ in $H$.

Homotope $\gamma^{\prime}$ in $H$ to eliminate the two points of intersection between the closure $\overline{\kappa^{\prime}}$ of $\kappa^{\prime}$ and the disc $D$.  Choose the homotopy so that $\gamma^{\prime}$ stays fixed outside a small regular neighborhood of $\overline{\kappa^{\prime}}$.  This reduces $|\gamma^{\prime}\cap D|$ by two and keeps all other points of $\gamma^{\prime}\cap D$ the same.  See Figure~\ref{homotopy} for an illustration.

\begin{figure}
\epsfig{file=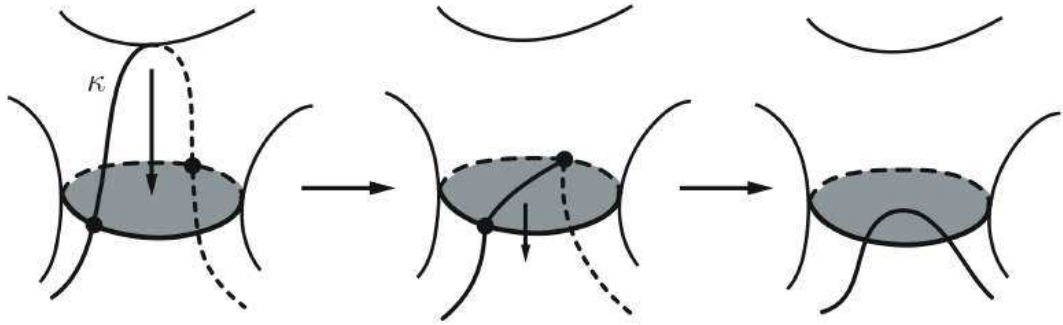, height=1.75in}
\caption{
We show a homotopy of $\gamma^{\prime}$ near the arc $\kappa$.  In the beginning, the arc $\kappa$ connects the two intersection points of $\gamma^{\prime}$ and the shaded disc $D$.  The homotopy reduces $|\gamma^{\prime}\cap D|$ by two and keeps $\gamma^{\prime}$ fixed outside a small neighborhood of $\overline{\kappa}$.
}
\label{homotopy}
\end{figure}

Consider the path $C$ in $\Gamma/\langle\gamma\rangle$.  Collapse the backtracking in $C$ to create a new path in $\Gamma/\langle\gamma\rangle$ with two fewer edges.  The new path $C$ represents a lift of the new path $\gamma^{\prime}$ to $\Gamma/\langle\gamma\rangle$ .

We now modify the disjoint collection of simple closed curves $B$ in $\hh$.  See Figure~\ref{recursionstep} for an illustration of this construction.  The covering map $\pi$ embeds $B$ into $\bd H$.  Consider the two points of $\gamma^{\prime}\cap D$ corresponding to the backtracking we chose above.  These two points correspond to two points of $\gamma\cap D$, and they lie in $\pi(B)$.  Consider the arc $\kappa$ of $\gamma-D$ that connects these two points in $\bd H$ and is homotopically trivial in $H \rel D$.  By construction, the arc $\kappa$ satisfies the property that its closure $\overline{\kappa}$ intersects $\pi(B)$ only at the endpoints of $\overline{\kappa}$.

\begin{figure}
\epsfig{file=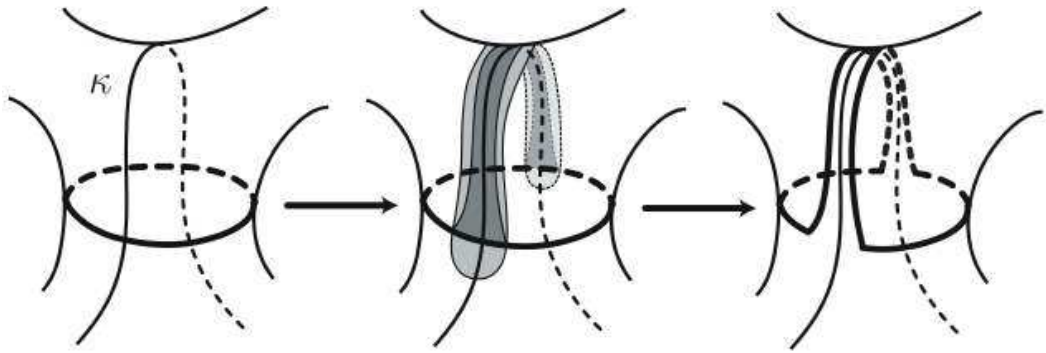, height=2in}
\caption{
The curve $\gamma$ intersects the thick curve $B$ on the left in a homotopically trivial arc $\kappa$.  We alter $B$ by considering a thin neighborhood $U$ around the arc $\kappa$, shown in light gray.  Consider a pair of arcs in $U$ running alongside $\kappa$.  Now cut and paste to form the new collection $B$.  The new collection is homologous to the old because the four arcs bound the rectangular region shown in darker gray.
}
\label{recursionstep}
\end{figure}

Let $U$ in $\bd H$ be a neighborhood of $\kappa$ chosen small enough so that (1) $U$ intersects $\pi(B)$ only in a small arc around each endpoint of $\overline{\kappa}$, and (2) the intersection $U\cap\gamma$ is a single arc component and a small neighborhood of $\kappa$ in $\gamma$.  Lift this picture back to $\hh$ with $\pi^{-1}$.  Replace the two arcs of $B\cap U$ in $U$ with a homologous pair of arcs that run parallel to $\kappa$. Choose the pair arcs so they do not intersect $\pi^{-1}(\gamma)$.

This gives us a new disjoint collection of simple closed curves $B$ in $\bd \hh$.  The collection $B$ satisfies the following three properties.  The new collection $B$ represents the same element of $H_{1}(\bd \hh, \Z)$ as before.  The map $\pi$ embeds $B$ into $H$.  The collection $\pi(B)$ intersects $\gamma$ in two fewer points than before.

Repeat the induction stage by returning to (\textbf{*}).

\textbf{**End induction stage.}

We now have a path $C$ in $\Gamma/\langle\gamma\rangle$ without backtrackings, a disjoint collection of simple closed curves $B$ in $\bd\hh$, and a curve $\gamma^{\prime}$ in $H$ homotopic to $\gamma$.  The edges of $C$ in $\Gamma/\langle\gamma\rangle$, the points of $\gamma\cap\pi(B)$ in $\hh$, and the points of $\gamma^{\prime}\cap D$ in $H$ are all in correspondence.

Lift $\gamma^{\prime}$ to a path $P^{\prime}$ in $\Gamma$.  The path $C$ is the projection of $P^{\prime}$ to $\Gamma/\langle\gamma\rangle$.  Since $C$ has no backtrackings, the path $P^{\prime}$ has no backtrackings.  By Lemma~\ref{geodesic}, the curve $\gamma^{\prime}$ minimizes intersection with $D$ over its homotopy class.  By Lemma~\ref{wedge}, the intersection number $|\gamma^{\prime}\cap D|$ equals $\girth(\gamma^{\prime})$.  Since $\gamma^{\prime}$ is homotopic to $\gamma$, the intersection number equals $\girth(\gamma)$.  Hence $C$ has $n$ edges.  Since the edges of $C$ are in correspondence with $\gamma\cap\pi(B)$, the curve $\gamma$ intersects $\pi(B)$ in $n$ points.

Since $B$ represents a nontrivial homology class in $\bd \hh$, at least one component of $B$ represents a nontrivial homology class.  Choose one of these components of $B$, and call it $B^{\prime}$.  Map $B^{\prime}$ to $H$ with $\pi$.  The curve $\pi(B^\prime)$ is a simple closed curve in $\bd H$ that bounds an embedded disc $D^\prime$ in $H$ and intersects $\gamma$ in no more than $n$ points.

The disc $D^{\prime}$ is essential because $B^{\prime}$ represents a nontrivial homology class in $\bd \hh$.  The disc $D^\prime$ intersects $\gamma$ in no more than $n$ points.  Since the girth of $\gamma$ is $n$, the disc $D^\prime$ intersects $\gamma$ in at least $n$ points.  The disc $D^{\prime}$ intersects $\gamma$ in exactly $n$ points, so we are done.
\end{proof}

\begin{figure}
\epsfig{file=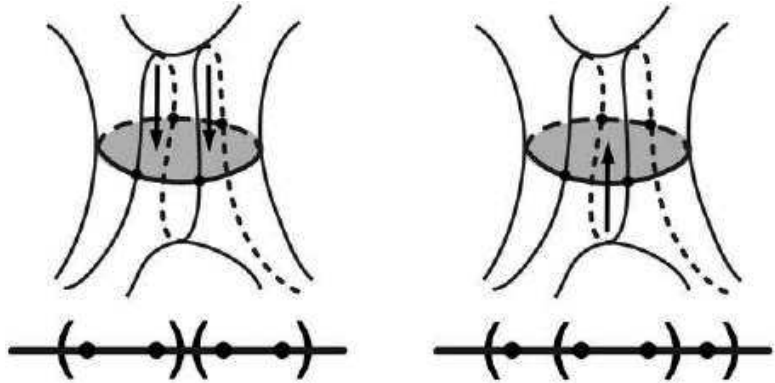, height=2in}
\caption{
The curve $\gamma$ intersects the shaded disc $D$ in four points.  The corresponding lift of $\gamma$ to $\Gamma$ has three backtrackings.  We illustrate two orders in which the backtrackings can be canceled, as in our second proof of Lemma~\ref{one embedded curve}.  The arrows depict the first cancellations.}
\label{parenthesize}
\end{figure}

The proof above constructs a geometrically girth-realizing disc $D^{\prime}$ from an algebraically girth-realizing disc $D$.  The construction differs depending on the order in which the backtrackings of $C$ were canceled.  In general, different sequences are possible.  Figure~\ref{parenthesize} shows a simple example with two possibilities.  We encourage the reader to construct the two resulting curve collections $B$.

\section{Girth of Heegaard Splittings}
\label{heegaardgirth}
A Heegaard splitting of a closed 3-manifold can be represented by a collection of disjoint simple closed curves on the boundary of a handlebody.  In the previous section, we examined the girth of simple closed curves on the boundary of a handlebody.  In this section we turn to girth for Heegaard splittings.  We define the algebraic and geometric girths of a Heegaard splitting.  We present Conjecture~\ref{ag girth conjecture} and prove that it is equivalent to the Girth Conjecture.  We also present some connections to the study of 3-manifolds with finite fundamental group.

Let $M$ be an orientable closed 3-manifold.  A \emph{Heegaard splitting} 
$H\cup H^\prime$ of $M$ is an embedding of a handlebody $H$ into $M$ so that $M-H$ is another handlebody $H^\prime$.  The intersection of the handlebodies in $M$ is the surface $\bd H$.  The order matters in our definition, so the splitting $H\cup H^\prime$ is different from $H^\prime\cup H$.  

Let $H\cup H^\prime$ be a Heegaard splitting of $M$.  The fundamental 
group $\pi_1(M)$ equals the quotient of $\pi_1(H)$ by some normal 
subgroup of relators $\rr$.
Each element of $\rr$ can be represented by an immersed 
curve in $\partial H$ that bounds a disc in $H^\prime$.
By the following lemma, the trivial element of $\rr$ can 
always be represented by an essential closed curve in $\bd H$.
This is why we require primitive curves $\gamma$ to be nontrivial in our 
statement of the Girth Conjecture.
The following lemma appears as Theorem~10.2 in \cite{P}.
We present a short proof here for clarity.

\begin{lemma}
\label{trivial relator}
Let $M$ be a closed orientable $3$-manifold.  Let $H\cup H^\prime$ 
be a Heegaard splitting of $M$ with genus strictly bigger than one.  
Then there is an essential closed curve in $\partial H$ that 
contracts in both $H$ and $H^\prime$.
\end{lemma}

\begin{proof}
Pick a basepoint $x$ in $\partial H$.  Choose an essential loop
$g$ (resp. $h$) based at $x$ that contracts in $H$ (resp. $H^\prime$).
Since $\genus(\partial H)>1$, we can choose $g$ and $h$ not to commute
in $\partial H$.  The commutator $[g,h]$ is an essential 
loop in $\partial H$ that contracts in both $H$ and $H^\prime$.
\end{proof}

Define the \emph{algebraic girth} of $H\cup H^\prime$ to be 
the minimum girth in $H$ over all nontrivial elements of $\rr$.  
When $\rr$ is trivial, we define the algebraic girth to be zero.  
Observe that $\rr$ is trivial if and only if $M$ is a connect 
sum of $g$ $S^2 \times S^1$'s, where $g$ is the genus of $H$.
The definition of algebraic girth is not symmetric because
it depends on a choice of one of the two handlebodies.  
It is not clear \emph{a priori} whether the algebraic girths of $H\cup H^{\prime}$ and $H^\prime \cup H$ are the same .

By the following lemma, higher genus Heegaard splittings of 3-manifolds with finite fundamental group have algebraic girth zero.

\begin{lemma}
\label{finite implies girth zero}
Let $M$ be a closed orientable 3-manifold.  Let $H\cup H^{\prime}$ be a Heegaard splitting of $M$.  If $M$ has finite fundamental group and $\genus(H)\geq 2$, then the algebraic girth of $H\cup H^{\prime}$ is zero.
\end{lemma}

\begin{proof}
Let $\dd$ be a disc system for $H$.  Let $\{ a, b, \ldots \}$ be the corresponding generating system for $\pi_{1}(H)$.  Consider the generator $a$.  Represent $a$ as a closed curve disjoint from the disc corresponding to the generator $b$.  Since $\pi_1(M)=\pi_1(H)/\rr$ is finite, there is an $n\geq 1$ so that $a^n$ is trivial in $\pi_1(M)$.  Since $a^n$ has girth zero in $H$, the splitting $H\cup H^\prime$ has algebraic girth zero.
\end{proof}

Let $\scc$ denote the collection of simple closed curves in $\partial H$ that bound a disc in $H^\prime$.  Define the \emph{geometric girth} of $H\cup H^\prime$ to be the minimum girth in $H$ over all essential elements of $\scc$.  By Lemma~\ref{one embedded curve}, the geometric girth equals the minimum intersection number of a pair of essential discs on either side of $\partial H$.  For example, a splitting is weakly reducible if and only if it has geometric girth zero.  Observe that geometric girth is symmetric with respect to $H$ and $H^{\prime}$.

The algebraic girth of a Heegaard splitting $H\cup H^\prime$ is less than or equal to its geometric girth.  We state this as a corollary of the following lemma.

\begin{lemma}
\label{girth zero}
Let $M$ be a closed orientable $3$-manifold.  Let $H\cup H^\prime$ 
be a Heegaard splitting of $M$.  If $H\cup H^\prime$ is reducible, then the algebraic girth of $H\cup H^\prime$ is zero.
\end{lemma}

Recall that a Heegaard splitting $H\cup H^{\prime}$ is reducible if there exists an essential  simple closed curve $\gamma$ in $\bd H$ that bounds a disc in both $H$ and $H^{\prime}$.

\begin{proof}
Since $H\cup H^{\prime}$ is reducible, choose a simple closed curve $\gamma$ in $\bd H$ that bounds a disc $D$ in $H$ and a disc $D^\prime$ in $H^\prime$.  Complete the disc $D^\prime$ in $H^\prime$ to a disc system $\dd^{\prime}$ in $H^\prime$.  Then the curves in $\bd \dd^{\prime}$ normally generate the group of relators $\rr$.  If all the curves in $\bd\dd^{\prime}$ are trivial in $\pi_1(H)$, then the group of relators $\rr$ is trivial; and the algebraic girth is zero by definition.

Otherwise, at least one curve in $\bd\dd^{\prime}$ is nontrivial in $\pi_1(H)$.  Consider such a curve.  This curve is disjoint from $D$, so it has girth zero in $H$.  It represents a nontrivial element of $\rr$, so the splitting has algebraic girth zero.
\end{proof}

We now prove a corollary of Lemma~\ref{girth zero}.

\begin{corollary}
\label{easy inequality}
Let $M$ be a closed orientable 3-manifold.  Let $H\cup H^\prime$ be a Heegaard splitting of $M$.  The algebraic girth of $H\cup H^\prime$ is less than or equal to its geometric girth.
\end{corollary}

\begin{proof}
Choose a simple closed curve $\gamma$ in $\bd H$ that realizes 
the geometric girth of $H\cup H^\prime$. Then $\gamma$ bounds a disc $D^\prime$ in $H^\prime$, and the girth of $\gamma$ in $H$ equals the geometric girth of $H\cup H^\prime$.

If $\gamma$ does not bound a disc in $H$, then $\gamma$ represents a nontrivial element of $\rr$.  Viewing $\gamma$ as an element of $\rr$, we conclude that the algebraic girth of $H\cup H^\prime$ is less than or
equal to the geometric girth.  If $\gamma$ does bound a disc $D$ in $H$, then $H\cup H^{\prime}$ is reducible.  Applying Lemma~\ref{girth zero}, we are done.
\end{proof}

The Girth Conjecture is equivalent to the assertion that the
reverse inequality is also true, except for a special case.  
We state the assertion here as a conjecture, and we prove its
equivalence below as Proposition~\ref{a and b}.

\begin{conjecture}
\label{ag girth conjecture}
Let $H\cup H^\prime$ be a Heegaard splitting of a closed orientable 
$3$-manifold.  If $H\cup H^\prime$ is not a genus two splitting of a genus two spherical manifold, then the algebraic and geometric girths of $H\cup H^\prime$ are equal.
\end{conjecture}

This is not necessary for our discussion, but the algebraic and geometric girths of genus two splittings of genus two spherical manifolds are not equal.  These splittings have algebraic girth zero and geometric girth two.  They have algebraic girth zero by Lemma~\ref{finite implies girth zero}.  They have geometric girth two for the following reasons.

Splittings of Seifert-fibered spaces have been classified in~\cite{MS}.  These splittings can be verified to have geometric girth at least two or less.  Genus two Heegaard splittings of geometric girth zero are all connect sums of Lens spaces, so their girth cannot be zero.  Finally, geometric girth one splittings of genus two or more do not exist by Lemmas~\ref{one embedded curve} and~\ref{geometric girth one}.

We prove Conjecture~\ref{ag girth conjecture} is equivalent to the Girth Conjecture at the end of this section.  For reference, we state their equivalence as a proposition here.

\begin{proposition}
\label{a and b}
Let $H\cup H^{\prime}$ be a Heegaard splitting of a closed orientable 3-manifold.  The Girth Conjecture is true for $H\cup H^{\prime}$ if and only if Conjecture~\ref{ag girth conjecture} is true for $H\cup H^{\prime}$.
\end{proposition}

It is often easier to use Conjecture~\ref{ag girth conjecture} in proofs than to cite the Girth Conjecture directly.  Here are three examples.

\begin{lemma}
\label{reducible}
Let $M$ be a closed orientable 3-manifold.  Let $H\cup H^\prime$ be a Heegaard splitting of $M$.  If $M$ is reducible, or more generally if $H\cup H^{\prime}$ is reducible, then the Girth Conjecture is true for $H\cup H^{\prime}$.
\end{lemma}

\begin{proof}
By~\cite{Ha}, the splitting $H\cup H^{\prime}$ is reducible if $M$ is reducible.  So consider a reducible Heegaard splitting $H\cup H^{\prime}$.

Since $H\cup H^{\prime}$ is reducible, it has geometric girth zero.  By Lemma~\ref{girth zero}, the splitting has algebraic girth zero.  Since the algebraic and geometric girths are equal, Conjecture~\ref{ag girth conjecture} is true for $H\cup H^{\prime}$.  By Proposition~\ref{a and b}, the Girth Conjecture is true for $H\cup H^{\prime}$.
\end{proof}

The Girth Conjecture for 3-manifolds with finite fundamental group is equivalent to the spherical space form conjecture.  The Geometrization Conjecture, then, implies the Girth Conjecture for Heegard splittings of 3-manifolds with finite fundamental group.  Conversely, the Girth Conjecture would provide another proof of the Geometrization Conjecture for 3-manifolds with finite fundamental group.

\begin{proposition}
\label{finite fundamental group}
Let $M$ be a closed orientable 3-manifold with finite fundamental group.  Let $H\cup H^\prime$ be an irreducible Heegaard splitting of $M$.  The 3-manifold $M$ is spherical if and only if the Girth Conjecture is true for $H\cup H^\prime$.
\end{proposition}

\begin{proof}
\textbf{S $\Rightarrow$ GC:}  
Let $M$ be a spherical 3-manifold.  By Lemma~\ref{reducible}, if $H\cup H^{\prime}$ is reducible, then the Girth Conjecture is true for $H\cup H^{\prime}$.  So assume the splitting $H\cup H^{\prime}$ is irreducible.

Irreducible Heegaard splittings of spherical 3-manifolds have been classified, and they are all genus two or less.  If $H$ has genus two, then $H\cup H^{\prime}$ satisfies Conjecture~\ref{ag girth conjecture}.  By Proposition~\ref{a and b}, it must also satisfy the Girth Conjecture.

So assume $H$ has genus one.  In this case, $H\cup H^{\prime}$ is a standard genus one splitting of a Lens space.  These have equal algebraic and geometric girths.  Consequently, they satisfy Conjecture~\ref{ag girth conjecture}, and hence the Girth Conjecture.

\textbf{GC $\Rightarrow$ S:}  
Let $H\cup H^{\prime}$ be a Heegaard splitting of $M$ that satisfies the Girth Conjecture.  By Proposition~\ref{a and b}, Conjecture~\ref{ag girth conjecture} is also true for $H\cup H^{\prime}$.

If the genus of $H$ is one, then $M$ is a Lens space and so is a spherical 3-manifold.  So assume that the genus of $H$ is at least two.  By Conjecture~\ref{ag girth conjecture}, either $M$ is a genus two spherical 3-manifold, or the algebraic and geometric girths of $H\cup H^\prime$ are equal.  In the former case, we are done.  So assume that the algebraic and geometric girths are equal.

By Lemma~\ref{finite implies girth zero}, the splitting $H\cup H^{\prime}$ has algebraic girth zero.  Therefore $H\cup H^{\prime}$ has geometric girth zero.  In other words, $H\cup H^\prime$ is weakly reducible.  By a theorem of Casson and Gordon \cite{CG}, the 3-manifold $M$ is either Haken or the splitting is 
reducible.  Since $M$ has finite fundamental group, the manifold cannot 
be Haken.  Thus $H\cup H^\prime$ is reducible.
\end{proof}

Recall Lemma~\ref{geometric girth one} which states that no simple closed curve on the boundary of a higher genus handlebody has girth one.  Equivalently, any Heegaard splitting of geometric girth one is the standard genus one splitting of the 3-sphere.  The same is true for algebraic girth.  Together these facts imply that the Girth Conjecture is true for all Heegaard splittings of algebraic girth one.  We state this as a proposition below.  For comparison, note that Conjecture~\ref{ag girth conjecture} and the Girth Conjecture are not known even for Heegaard splittings of algebraic girth zero.

\begin{proposition}
\label{algebraic girth one}
Let $H\cup H^\prime$ be a Heegaard splitting of a closed 3-manifold.  If $H\cup H^\prime$ has algebraic girth one, then the splitting is the standard genus one splitting of the 3-sphere.  In particular, the Girth Conjecture holds for $H\cup H^\prime$.
\end{proposition}

Note that Proposition~\ref{algebraic girth one} also follows from Lemma~\ref{group girth one}.

\begin{proof}
Let $H\cup H^\prime$ be a Heegaard splitting of a closed orientable 3-manifold.  Let $\rr$ denote the group of relators, so that $\pi_1(M)=\pi_1(H)/\rr$.  Assume that $\rr$ has a relator of girth one.

It suffices to show that either $H\cup H^\prime$ is a genus one splitting of the 3-sphere, or $H\cup H^\prime$ has algebraic girth zero.  For if $H\cup H^{\prime}$ were the genus one splitting of the 3-sphere, then it would have geometric girth one.  This would imply Conjecture~\ref{ag girth conjecture} for $H\cup H^{\prime}$.  By Proposition~\ref{a and b}, this would also imply the Girth Conjecture for $H\cup H^{\prime}$.

Let $r$ be an element of $\rr$ of girth one.  Let $\kappa$ be a curve in $H$ representing $r$.  Then $\kappa$ can be made to intersect some essential disc $D$ in $H$ exactly once.  The disc $D$ must be non-separating.

Let $n$ be the genus of $H$.  If $n$ equals one, then $\pi_1(M)$ is trivial and $H\cup H^\prime$ represents the 3-sphere.  So assume that $n$ is greater than one.  By Lemma~\ref{genus two}, the genus of $H$ cannot equal two, so we have $n\geq 3$.

Extend $D$ to a disc system $\dd$ for $H$.  Let $\{a_{1}, a_{2}, \ldots a_{n}\}$ denote the corresponding generating system for $\pi_{1}(H)$, with $a_{1}$ corresponding to $D$.  Homotope $\kappa$ to minimize intersection with $\dd$.  By Corollary~\ref{subcollection}, the curve $\kappa$ intersects $D$ exactly once.  Let $r=a_{1}\cdot w(a_{2}, a_{3}, \ldots a_{n})$ denote the corresponding word in the generators, where $a_{1}$ does not occur in the word $w$.  Then $r$ is cyclically reduced.

Choose a disc system $\dd^{\prime}$ for $H^\prime$.  Then $\bd\dd^{\prime}$ gives a collection of $n$ relators $g_i$ that normally generate $\rr$.
Multiplying by conjugates of $r$, we can eliminate every occurrence of $a_{1}$ from the $g_{i}$.  This gives $n$ new relators $g_i$ that, together with $r$,
normally generate $\rr$.  If one of the $g_{i}$, say $g_{1}$, is nontrivial in $\pi_{1}(H)$, then $H\cup H^{\prime}$ has algebraic girth zero.  This follows because $g_{1}$ would be a nontrivial relator represented by a curve disjoint from $D$.

Assume that all of the $g_{i}$ are trivial in $\pi_{1}(H)$.  Then $\rr$ is normally generated by the relator $r$.  We prove that $\pi_{1}(M)$ equals the free group generated by $\{ a_{2}, a_{3}, \ldots a_{n}\}$.  Let $\Phi$ be the natural map $\Phi:\langle a_{2}, a_{3}, \ldots a_{n} \rangle\rightarrow\pi_{1}(H)/\rr=\langle a_{1}, a_{2}, \ldots a_{n} \rangle/ (r)$.  The map $\Phi$ is surjective because $a_{1}=w(a_{2}, a_{3}, \ldots a_{n})^{-1}$ in $\pi_{1}(H)/\rr$.  By the Freiheitssatz~\cite[Thm. 4.10]{MKS}, every nontrivial relator in the normal subgroup generated by $r$ must involve $a_{1}$.  It follows that $\Phi$ is injective.

Consequently, $\pi_{1}(M)$ is isomorphic to a free group on $n-1$ letters.  Since $n\geq 3$, the fundamental group $\pi_{1}(M)$ is equal to a nontrivial free product $A\star B$.  By Kneser's Conjecture~\cite[Theorem II.A.3]{St2}, the 3-manifold $M$ is reducible.  By~\cite{Ha}, this implies the Heegaard splitting $H\cup H^{\prime}$ is reducible.  Finally, applying Lemma~\ref{girth zero}, the algebraic girth of $H\cup H^{\prime}$ is zero.
\end{proof}

We now prove Proposition~\ref{a and b}.

\begin{proof}[Proof of Proposition~\ref{a and b}]
Let $H\cup H^\prime$ be a Heegaard splitting of a closed
orientable 3-manifold $M$.

\textbf{\ref{ag girth conjecture} $\Rightarrow$ GC:}  
Assume Conjecture~\ref{ag girth conjecture} for $H\cup H^{\prime}$.
Let $\gamma$ be a closed curve in $\bd H$ that is primitive 
in $\pi_1(H)$ and contracts in $H^\prime$.  It suffices to construct a simple closed curve $\gamma^{\prime}$ in $\bd H$ that bounds a disc in $H^{\prime}$ and satisfies $\girth(\gamma^{\prime}) \leq \girth(\gamma)$.

Since $\gamma$ is primitive, the curve $\gamma$ represents a nontrivial element of $\rr$.  In particular, the algebraic girth of $H\cup H^\prime$
is less than or equal to $\girth(\gamma)$.
Choose a simple closed curve $\gamma^\prime$ in $\bd H$ that 
realizes the geometric girth of $H\cup H^\prime$.  
We will show that $\girth(\gamma^{\prime}) \leq \girth(\gamma)$.  By Conjecture~\ref{ag girth conjecture}, either $H\cup H^{\prime}$ is a genus two splitting of a genus two spherical manifold, or the algebraic and geometric girths are equal.

Assume first that the algebraic and geometric girths are equal.  Then $\girth(\gamma^{\prime})$ equals the algebraic girth of $H\cup H^\prime$.  Because the algebraic girth of $H\cup H^\prime$ is less than or equal to $\girth(\gamma)$, we are done in this case.

Assume next that $H\cup H^\prime$ is a genus two splitting of 
a genus two spherical manifold.  Then $\pi_1(M)$ is not a 
cyclic group.  Moreover, the geometric girth of $H\cup H^\prime$
is exactly two.  This is by our remarks following the statement of Conjecture~\ref{ag girth conjecture}.  Hence $\girth(\gamma^\prime)=2$.  It suffices to show that $\gamma$ cannot have girth zero or one.

In a genus two handlebody, the only primitive curves of girth zero 
correspond to group elements of the form $a$, where $\langle a, b \rangle$ is a 
generating system for $\pi_1(H)$.  The curve $\gamma$ cannot be of this form because this would imply $\pi_1(M)$ is cyclic.  Thus $\gamma$ does not have zero girth.  Lastly, the curve $\gamma$ cannot have girth one by Lemma~\ref{genus two}.  So we are done with this implication.

\textbf{GC $\Rightarrow$~\ref{ag girth conjecture}:}
Now assume the Girth Conjecture for $H\cup H^\prime$.  It suffices to show that either $H\cup H^\prime$ is a genus two splitting of a genus two spherical manifold, or the algebraic and geometric girths of $H\cup H^\prime$ are equal.

If the genus of $H$ is one, then $H\cup H^\prime$ is a standard Heegaard splitting of a Lens space.  It is straightforward in this case that the algebraic and geometric girths are equal.  So assume that the genus of $H$ is at least two.

We have $\pi_1(M)=\pi_1(H)/\rr$, where $\rr$ is the group of relators.  If $H\cup H^\prime$ is a reducible Heegaard splitting, then the geometric girth of 
$H\cup H^\prime$ is zero.  In this case, the algebraic and geometric 
girths must be equal.

Assume then that $H\cup H^\prime$ is an irreducible Heegaard splitting. By~\cite{Ha}, the 3-manifold $M$ is irreducible, thus $\rr$ is nontrivial.  Let $\gamma$ be a closed curve in $\bd H$ that realizes the algebraic girth of the splitting $H\cup H^\prime$.  Then $\gamma$ is a nontrivial element of $\rr$, and $\gamma$ bounds a disc in $H^\prime$.  By Corollary~\ref{easy inequality}, the algebraic girth of $H\cup H^{\prime}$ is less than or equal to the geometric girth.    To show the algebraic and geometric girths are equal, it suffices to show the geometric girth is less than or equal to $\girth(\gamma)$.

Assume temporarily that $\gamma$ is primitive in $\pi_1(H)$.  Apply the Girth
Conjecture to $\gamma$.  Then there exists a simple closed curve $\gamma^\prime$ in $\bd H$ that bounds a disc in $H^\prime$ and satisfies $\girth(\gamma^{\prime})\leq \girth(\gamma)$.  Since the geometric girth of $H\cup H^\prime$ is less than or equal to $\girth(\gamma^\prime)$, we know the geometric girth is less than or equal to $\girth(\gamma)$.  So we are done in this case.

Assume from now on that $\gamma$ is not primitive in $\pi_1(H)$.  Then 
$\gamma$ is homotopic to $\alpha^n$ for some primitive $\alpha$ in $\pi_1(H)$.  By Lemma~\ref{power},
$\girth(\alpha)\leq\girth(\gamma)$.  If $\alpha$ is also in $\rr$, we 
can apply the Girth Conjecture to $\alpha$ as in the previous 
paragraph to be done.  So assume that $\alpha$ does not lie in $\rr$.
Then $\pi_1(M)$ has torsion.  Since $M$ is irreducible, we must 
have that $\pi_1(M)$ is finite~\cite[Chapter 9]{He1}.

We now treat the cases of $\genus(H)=2$ and $\genus(H)>2$, separately.
Assume first that $\genus(H)=2$.  By the work of Thurston, genus two manifolds are known to satisfy the Geometrization Conjecture~\cite{Thu1}. Since $\pi_{1}(M)$ is finite, the manifold $M$ is a spherical 3-manifold. Since $H\cup H^\prime$ is irreducible, the splitting must be a genus two splitting of a genus two spherical manifold.

Assume finally that $\genus(H)\geq 3$.  Represent $\pi_1(H)$ as a free 
group $\langle a, b, c, \cdots\rangle$ on finitely many generators
by choosing a disc system for $H$ and letting each disc correspond to
a generator.  Since $\pi_1(M)=\pi_1(H)/\rr$ is finite, there exist 
$m,n\geq 1$ so that both $a^m$ and $b^n$ are trivial in $\pi_1(M)$.  
Then $a^mb^n$ is primitive element in $\pi_1(H)$ and trivial in $\pi_{1}(M)$.  Since $a^mb^n$ can be represented by a curve in $H$ disjoint from 
the disc corresponding to $c$, the element $a^m b^n$ has girth zero.

Represent $a^mb^n$ by a curve in $\bd H$ that contracts in $H^\prime$.  Applying the Girth Conjecture to this curve, there exists a simple closed curve in $\bd H$ of girth zero in $H$ that bounds a disc in $H^\prime$.  In other words, the splitting $H\cup H^\prime$ has geometric girth zero.  By Corollary~\ref{easy inequality}, $H\cup H^\prime$ must also have algebraic girth zero.  The algebraic and geometric girths of $H\cup H^{\prime}$ are equal, so we are done.
\end{proof}

\section{Discs and $\dd$-Interfaces}
\label{subordination}
Let $H$ be a handlebody.  Let $\dd$ be a disc system for $H$.  Let $D$ be an essential disc in $H$.  The disc $D$ may be very complicated with respect to $\dd$.  In this section we prove Lemma~\ref{subordination-lemma}.  Lemma~\ref{subordination-lemma} lets us construct from $D$ either a disc disjoint from $\dd$ or two ``half-discs'' in $H-\dd$.  In both cases, the constructed objects preserve information about $D$.  Lemma~\ref{subordination-lemma} will be important for later sections.

Let $\dd$ be either a disc system for a handlebody $H$ or a generating system for a free group $F$.  Suppose $\dd$ has $g$ elements.  We use $\dpm$ to denote the set of $2g$ elements consisting of elements of $\dd$ with an orientation.  In the handlebody case, elements of $\dpm$ correspond to the $2g$ sides of discs in $\dd$.  We will usually refer to elements of $\dpm$ simply as discs.  In the free group case, the elements of $\dpm$ correspond to the generators of $F$ and their inverses.  The elements of $\dpm$ generate $F$ as a semigroup.

Define a $\dd$-\emph{interface} to be a partition of $\dpm$ into three subsets labeled black, white, and gray.  We denote these subsets $\db$, $\dw$, and $\dg$.  We require the black and white subsets to be nonempty.  This parallels the definition of interface in Section~\ref{trees}.  We use a single letter, bold and underlined, to denote a $\dd$-interface.  If we want to reference the black and white subsets of $\dpm$ explicitly, we use the notation $\di{\db}{\dw}$.  A \emph{$\dd$-disc} is a $\dd$-interface with $\dg$ empty.  A \emph{$\dd$-half-disc} is a $\dd$-interface with $\dg$ consisting of one element.  The connection to geometric discs and half-discs will become clear during the proof of Lemma~\ref{subordination-lemma}.

Let $H$ be a handlebody of genus $g$.  Let $\dd$ be a disc system for $H$.  Let $\dint{a}$ be a $\dd$-interface.  The $\dd$-interface $\dint{a}$ gives rise to an interface for $\hh$ in the following way.

Consider the tree $\Gamma$ associated to $\dd$.  The ends of $\Gamma$ are canonically in one-to-one correspondence with the ends of $\hh$.  There is an equivariant homotopy $\Phi$ from $\hh$ to $\Gamma$ such that the preimage of the midpoint of any edge equals a lift of the corresponding disc in $\dd$.  The homotopy $\Phi$ identifies the ends of $\hh$ with the ends of $\Gamma$.  

View $\Gamma$ as the Cayley graph of $\pi_{1}(H)$ with generating system $\dd$.  Choose a vertex $v$ of $\Gamma$.  The $\dd$-interface $\dint{a}$ gives colors to the $2g$ edges leaving $v$.  Delete $v$ from $\Gamma$ to form $2g$ connected subsets of $\Gamma$.  Each connected component inherits a color coming from $\dint{a}$.  These components partition the ends of $\Gamma$ and define an interface for $\Gamma$.  This carries over to an interface $\alpha$ for $\hh$.  The interface $\alpha$ is well-defined up to the choice of vertex $v$.  We say $\alpha$ is a \emph{lift} of $\dint{a}$ to $\hh$.  This parallels the notion of lifting a disc in $H$ to an interface for $\hh$.
       
Let $H$ be a handlebody, and let $\hh$ be the universal cover.  Recall from Section~\ref{trees} what it means for a proper curve to cross an interface.  Let $\alpha$ and $\beta$ be two interfaces for $\hh$.  The interface $\alpha$ is \emph{subordinate} to $\beta$ if, for every proper curve $\gamma$, the curve $\gamma$ crosses $\beta$ if $\gamma$ crosses $\alpha$.  In other words, the support of $\alpha$ is contained in the support of $\beta$.  Two interfaces $\alpha$ and $\beta$ are \emph{independent} if no proper curve crosses both.  In other words, they are independent if their supports are disjoint. 

Let $\dd$ be a disc system for $H$.  Let $D$ be a disc in $H$.  Recall what it means to lift a disc to an interface.  A $\dd$-interface $\dint{a}$ is \emph{subordinate} to $D$ if some lift of $\dint{a}$ to $\hh$ is subordinate to a lift of $D$.   Two $\dd$-interfaces $\dint{a}$ and $\dint{b}$ are \emph{independently subordinate}
to $D$ if they can be lifted to independent interfaces that are subordinate to a common lift of $D$.

We can now state and prove Lemma~\ref{subordination-lemma}.

\begin{lemma}
\label{subordination-lemma}
Let $H$ be a handlebody.  Let $\dd$ be a disc system for $H$.  Let $D$ be an essential embedded disc in $H$.  There is either a $\dd$-disc subordinate to $D$ or two $\dd$-half-discs independently subordinate to $D$.
\end{lemma}

\begin{proof}
Let $g$ be the genus of $H$.  Let $\hh$ be the universal cover of $H$.  Let $\Gamma$ be the infinite tree that is homotopic to $\hh$, as constructed above.  The ends of $\Gamma$ are in correspondence with the ends of $\hh$, up to a choice of deck transformation.  Choose such an identification.  Lift $D$ to an interface for $\hh$.  This defines an interface $\alpha$ for $\Gamma$, dividing the ends of $\Gamma$ into  black and white subsets.

Assign the labels black and white to the unoriented edges of $\Gamma$ as follows.  An unoriented edge can be assigned either one, two, or no labels.  Let $(x,y)$ be an oriented edge of $\Gamma$ with origin $x$ and terminal vertex $y$.  Define a \emph{ray} from $x$ to $y$ to be an infinite geodesic path whose initial edge has origin $x$ and terminus $y$.  Every ray in $\Gamma$ determines a unique end of $\Gamma$.

Label the unoriented edge corresponding to $(x,y)$ black (resp. white) if every ray from $x$ to $y$ determines a black (resp. white) end.  If an edge is labeled neither black nor white, label the edge gray.  Since $D$ lifts to a compact
set in $\hh$, the set of gray edges is finite.  Also, the set of gray edges is connected by an easy argument that follows from the definitions.

There are two cases.  Assume first that the set of gray edges is nonempty.  Consider the subgraph of gray edges.  Choose two outermost vertices of this subgraph.  Call the vertices $v$ and $w$.  The $2g$ edges leaving
$v$ are in correspondence with $\dpm$.  Since $v$ is outermost, there is exactly one gray edge leaving $v$.  The other $2g-1$ edges are not all black
and not all white, for otherwise the gray edge would be black or white.  This defines a $\dd$-half-disc that lifts to an interface based at the vertex $v$.  By construction, this interface is subordinate to $\alpha$.  Similarly, we have another $\dd$-half-disc that lifts to an interface based at $w$ and is subordinate to $\alpha$.  Since $v$ and $w$ are outermost vertices of a connected gray subgraph, the two lifts of $\dd$-half-discs are independent.  Thus they are independently subordinate to $D$.

Assume next that there are no gray edges.  Then every edge is colored either black, white, or both.  There are two subcases.  Assume first that some edge $(x,y)$ is both black and white.  Then all rays from $x$ to $y$ define black ends, say, and all rays from $y$ to $x$ define white ends.  Consider the vertex $x$.  Define a $\dd$-disc by setting the oriented edge corresponding to $(x,y)$ black, and the other $2g-1$ edges leaving $x$ white.  This lifts to an interface for $\Gamma$ at $x$ that is equal to $\alpha$, and hence subordinate to $\alpha$.

In the final subcase, each edge has a single label black or white.  Choose a vertex $x$ that touches both a black edge $(x,y)$ and a white edge $(x,z)$.  Then all the rays from $x$ to $y$ are black, and all the rays from $x$ to $z$ are white.  In fact, for any vertex $y$ adjacent to $x$, the rays from $x$ to $y$ are either all black or all white.  Define a $\dd$-disc based at $x$ accordingly.  This defines an interface for $\Gamma$ at $x$ that equals $\alpha$.  This interface is of course subordinate to $\alpha$ and is a lift of a $\dd$-disc.  
\end{proof}

\section{Computing the Girth of a Conjugacy Class}
\label{curvealgorithm}
\label{computecurve}
In this section we present an effective algorithm to compute the girth of a conjugacy class in a free group.  This gives an algorithm to compute the girth of a closed curve in a handlebody $H$.  To compute the girth of a curve, simply pass to the level of the fundamental group $\pi_{1}(H)$.  Using the results of this section, we also show that free groups of rank bigger than one do not contain elements of girth one.  We prove this at the end.

In the previous section we introduced the notion of a $\dd$-interface.  Let $H$ be a handlebody.  Let $\dd$ be a disc system for $H$.  Let $\gamma$ be a closed curve in $H$, and let $\dint{a}$ be a $\dd$-interface.  We define the intersection number $\gamma\wedge\dint{a}$ as follows.  Lift $\dint{a}$ to an interface $\alpha$ for $\hh$, as in Section~\ref{subordination}.  This lift is well-defined up to deck transformations.  Consider all lifts of $\gamma$ to $\hh$.  Define the intersection number $\gamma\wedge\dint{a}$ to be the number of lifts of $\gamma$ that cross the interface $\alpha$.  By equivariance, this quantity does not depend on the particular lift of $\dint{a}$.  The following lemma says how we can compute the intersection number.

\begin{lemma}
\label{intersection handlebody}
Let $H$ be a handlebody.  Let $\dd$ be a disc system for $H$.  Let $\gamma$ be a closed curve in $H$.  Let $\dint{a}$ be a $\dd$-interface.  Suppose $\gamma$ minimizes intersection with $\dd$ in its homotopy class.  The disc system $\dd$ divides $\gamma$ into a collection of arcs connecting elements of $\dpm$.  The intersection number $\gamma\wedge\dint{a}$ equals the number of arcs of $\gamma$ connecting a black element of $\dpm$ to a white element.
\end{lemma}

\begin{proof}
Let $g$ be the genus of $H$.  Form the directed graph $\Gamma_{\dd}$ consisting of a single vertex and $g$ directed edges in correspondence with $\dd$.  There exists a homotopy $\Phi:H\rightarrow\Gamma_{\dd}$.  Choose $\Phi$ so that the preimage under $\Phi$ of the midpoints of the edges in $\Gamma_{\dd}$ equals $\dd$.

Consider the universal covers of $H$ and $\Gamma_{\dd}$.  Denote them $\hh$ and $\Gamma$.  The homotopy $\Phi$ lifts to a homotopy from $\hh$ to $\Gamma$.  The graph $\Gamma$ is the tree associated to $\dd$.  We can also view $\Gamma$ as the Cayley graph of $\pi_{1}(H)$ with respect to the generating system corresponding to $\dd$.  Lift $\gamma$ to a path in $\Gamma$.  By Lemma~\ref{geodesic}, the curve $\gamma$ lifts to a path without backtrackings.

Choose a vertex $v$ of $\Gamma$ and lift the $\dd$-interface $\dint{a}$ to an interface $\alpha$ for $\Gamma$.  The $2g$ edges that touch $v$ are in correspondence with $\dpm$.  Each of these $2g$ edges inherits a color black, white, or gray from $\dint{a}$.  The intersection number $\gamma\wedge\dint{a}$ equals the number of lifts of $\gamma$ to $\Gamma$ that cross $\alpha$.  Any geodesic path in $\Gamma$ that crosses $\alpha$ must pass through a black edge touching $v$, the vertex $v$, and a white edge touching $v$.  

Consider all lifts of $\gamma$ to $\Gamma$ that pass through $v$.  These lifts are in correspondence with the arcs $\gamma-D$ in $H$.  The lifts crossing $\alpha$ correspond to the arcs connecting a black element of $\dpm$ to a white element, so we are done.
\end{proof}

Let $F$ be a free group.  Let $\dd$ be a generating system for $F$.  Let $w$ be an element of $F$, and let $\dint{a}$ be a $\dd$-interface.  Then $\dint{a}$ assigns a color of black, white, or gray to each of the semigroup generators in $\dpm$.  In parallel with the above, we define the intersection number $w\wedge\dint{a}$.  Form a handlebody $H$ with disc system $\dd$.  Choose a basepoint for $H$.  The group $\pi_{1}(H)$ is naturally in correspondence with $F$.  Choose a closed loop $\gamma$ in $H$ that gives rise to $w$.  View $\dint{a}$ as a $\dd$-interface for the handlebody $H$.  Define $w\wedge\dint{a}$ as $\gamma\wedge\dint{a}$.  This definition does not depend on the choice of $H$, the choice of closed curve, or the choice of basepoint.  The following lemma is simply Lemma~\ref{intersection handlebody} phrased in terms of free groups.

\begin{lemma}
\label{intersection free group}
Let $F$ be a free group.  Let $\dd$ be a generating system for $F$.  Let $w$ be a word in the semigroup generators $\dpm$ of $F$.  Let $\dint{a}$ be a $\dd$-interface.  Suppose that $w$ is cyclically reduced.  The intersection number $w\wedge\dint{a}$ equals the number of adjacent pairs $ab$ in $w$ such that $a^{-1}$ is black and $b$ is white, or $a^{-1}$ is white and $b$ is black.
\end{lemma}

\begin{proof}
Realize $F$ as the fundamental group of a handlebody $H$ with disc system $\dd$, as above.  Consider any closed curve $\gamma$ in $H$ that is transverse to the discs in $\dd$.  By looking at the sequence of intersections with elements of $\dd$, the curve $\gamma$ gives rise to a word in the elements of $\dpm$.  In other words, $\gamma$ gives rise to a word in the generators $\dd$ and their inverses.

Let $\Gamma$ denote the tree associated to $\dd$.  Choose a closed curve $\gamma$ in $H$ that gives rise to the word $w$.  Since $w$ is cyclically reduced, the curve $\gamma$ has no backtrackings when lifted to $\Gamma$.  By Lemma~\ref{geodesic}, the curve $\gamma$ minimizes intersection with $\dd$ in its homotopy class.

The quantity $w\wedge\dint{a}$ equals $\gamma\wedge\dint{a}$ by definition.  By Lemma~\ref{intersection handlebody}, the quantity $\gamma\wedge\dint{a}$ equals the number of arcs of $\gamma-D$ connecting a black element of $\dpm$ to a white element.   An arc of $\gamma-D$ corresponds to an adjacent pair of letters in the word $w$.  The arcs of $\gamma-D$ that connect a black element of $\dpm$ to a white element correspond to the adjacent pairs of letters $ab$ such that $a^{-1}$ is black, resp. white, and $b$ is white, resp. black.
\end{proof}

Let $w$ be an element of the free group $F$.  Let $\dd$ be a disc system for $F$.  Define the \emph{complexity} of $\dd$ with respect to $w$ to be $w\wedge\dd$.  By Lemma~\ref{cyclic}, the complexity of $\dd$ equals the length of $w$ written as a cyclically reduced word in the generators corresponding to $\dd$.

Let $\dint{c}$ be a $\dd$-half-disc $\di{\db}{\dw}$.  Let $\dg$ the corresponding gray disc.  We say $\dint{c}$ is a \emph{shortcut} $\dd$-half-disc for $\dd$ with respect to $w$ if $w\hi\di{\db}{\dw}$ is strictly less than both $\gamma\hi\di{\db}{\dg}$ and $\gamma\hi\di{\dw}{\dg}$.  Note that $\gamma\hi\di{\db}{\dg}$, for instance, is the $\dd$-interface obtained by taking $\dd_{B}$ as the set of black discs and $\dd_{G}$ as the set of white discs.

By the following lemma, if $\dd$ has a shortcut $\dd$-half-disc, we can create a generating system of smaller complexity.  Thus we say $\dd$ is \emph{locally minimal} with respect to $w$ if $\dd$ has no shortcut $\dd$-half-discs.  Locally minimal disc systems are useful for calculating girth by Lemma~\ref{locally-minimal}.

\begin{lemma}
\label{curve shortcut}
Let $F$ be a free group.  Let $\dd$ be a generating system for $F$.  Let $w$ be an element of $F$.  If there is a shortcut $\dd$-half-disc with respect to $w$, then there is a generating system of smaller complexity with respect to $w$.
\end{lemma}

\begin{proof}
The set $\dpm$ generates $F$ as a semigroup.  Suppose $\di{\db}{\dw}$ is a shortcut $\dd$-half-disc with respect to $w$, with gray element $\dg$.   Then $w\hi\di{\db}{\dw}<\min(w\hi\di{\db}{\dg},w\hi\di{\dw}{\dg})$.  The inverse of the element corresponding to $\dg$ in $\dpm$ is either black or white.  Assume the inverse is black, say, so that it lies in $\db$.  We now work geometrically.

Realize the free group $F$ and generating system $\dd$ as the fundamental group of a handlebody $H$ with disc system labeled by $\dd$.  Choose a closed curve $\gamma$ that gives rise to $w$.  Then $\gamma\hi\di{\db}{\dw}<\min(\gamma\hi\di{\db}{\dg},\gamma\hi\di{\dw}{\dg})$.

Let $D_{G}$ denote the geometric disc in $\dd$ that corresponds to $\dg$.  Realize the $\dd$-disc $\di{\db}{\dg\cup\dw}$ as a new geometric disc $D_{N}$ in $H-\dd$.  Replace $D_{G}$ with $D_{N}$ to create a new disc system $\dd^\prime$.  The collection $\dd^\prime$ is a disc system because, by construction, $D_{N}$ is non-separating in $H$:  the collection $\db$ contains the element of $\dpm$ paired with $\dg$ in $H$.

The complexity of $\dd^{\prime}$ is the sum of $\gamma\hi D$ over all discs $D$ in $\dd$.  The discs in $\dd^{\prime}$ are the same as those in $\dd$, except that $D_{N}$ replaced $D_{G}$.  These two discs satisfy $\gamma\hi D_{N}=\gamma\hi\di{\db}{\dg}+\gamma\hi\di{\db}{\dw}<\gamma\hi\di{\db}{\dg}+\gamma\hi\di{\dw}{\dg}=\gamma\hi D_{G}$.  The middle inequality follows from the shortcut inequality $\gamma\hi\di{\db}{\dw}<\gamma\hi\di{\dw}{\dg}$.  Since $\gamma\hi D_{N}<\gamma\hi D_{G}$, the disc system $\dd^\prime$ has smaller complexity than $\dd$.
\end{proof}

The above proof is constructive.  It shows how to construct a disc system of smaller complexity if a shortcut $\dd$-half-disc exists.  The construction is geometric.  For algorithmic purposes we restate the construction algebraically below.

Let $F$ be a free group with generating system $\dd$.  The elements of $\dpm$ generate $F$ as a semigroup.  Let $w$ be a word in the generators $\dpm$.  Let $\di{\db}{\dw}$ be a shortcut $\dd$-half-disc with respect to $w$, with gray element $\dg$.  Cyclically reduce $w$.  By Lemma~\ref{cyclic}, the complexity of $\dd$ is the length of $w$.

We will create a new generating system $\dd^{\prime}$ by constructing an automorphism $\Phi$ of $F$.  Automorphisms preserve girth.  By Lemma~\ref{curve shortcut}, the word $\Phi(w)$, when cyclically reduced, will have length shorter than $w$ with respect to the generating system $\dd$.  By taking $\dd^{\prime}$ to be $\Phi^{-1}(\dd)$, we will obtain a generating system of smaller complexity with respect to $w$.

Denote as $c$ the group element in $\dpm$ corresponding to $\dg$.  Switch the black and white subsets of $\di{\db}{\dw}$, if necessary, so that $c^{-1}$ is black.  This mirrors the proof of Lemma~\ref{curve shortcut}.  If we translate the proof of Lemma~\ref{curve shortcut} into algebraic terms, we get the automorphism defined in Table~\ref{automorphism}.

\begin{table}
\begin{center}
\begin{tabular}{c|c|c|}
\multicolumn{3}{r}{black \ \ white} \\ \cline{2-3}
black & $c^{-1}ac$ & \ $c^{-1}a$ \ \\ \cline{2-3}
white & $ac$ & $a$ \\ \cline{2-3 }
\end{tabular}
\end{center}
\caption{This table gives instructions on where $\Phi$ should map each generator $a$ of the generating system $\dd$ for $F$.  The row corresponds to whether $a$ is black or white as an element of $\dpm$.  The column corresponds to whether $a^{-1}$ is black or white as an element of $\dpm$.}
\label{automorphism}
\end{table}

\begin{lemma}
\label{locally-minimal}
Let $F$ be a free group.  Let $\dd$ be a generating system for $F$.  Let $w$ be an element of $F$.  If $\dd$ is locally minimal with respect to $w$, then the girth of $w$ equals the minimum of $w\hi\dint{a}$ over all $\dd$-discs $\dint{a}$.
\end{lemma}

\begin{proof}
Form a handlebody $H$ with disc system labeled by $\dd$.  Choose a basepoint for $H$.  The fundamental group $\pi_{1}(H)$ is in natural correspondence with $F$.  Let $\gamma$ be a closed loop in $H$ that gives rise to $w$.  Any $\dd$-disc for $F$ transfers to a $\dd$-disc for $H$.  If $\dint{a}$ is a $\dd$-disc for $F$, then $w\hi\dint{a}$ equals $\gamma\hi\dint{a}$, with the second expression taken in $H$.  By definition, the girth of $w$ is the minimum of $\gamma\hi D$ over all essential discs $D$ in $H$.  Let $n$ be the girth of $\gamma$.  Let $D$ be a girth-realizing disc for $\gamma$, so that $\gamma\hi D=n$.  Apply Lemma~\ref{subordination-lemma} to the disc $D$ in $H$ with disc system $\dd$.  There are two cases.

\textbf{Case 1:  There is a $\dd$-disc subordinate to $D$.}
Let $\dint{a}$ be a $\dd$-disc subordinate to $D$.  Since $\dint{a}$ 
is subordinate to $D$, we have $\gamma\hi\dint{a}\leq\gamma\hi D$.
Thus $\gamma\hi\dint{a}\leq n$.  Proceed with $\dint{a}$ to the conclusion 
below.

\textbf{Case 2:  There are two $\dd$-half-discs independently 
subordinate to $D$.}  Since there are two $\dd$-half-discs independently subordinate to $D$, one of the $\dd$-half-discs, call it $\dint{c}$, satisfies 
$\gamma\hi\dint{c}\leq\lfloor\frac{n}{2}\rfloor$.  Let $\di{\db}{\dw}$ denote $\dint{c}$, with gray disc $\dg$.  Then $\gamma\hi\di{\db}{\dw}\leq\lfloor\frac{n}{2}\rfloor$.  Since $\dd$ is locally minimal, the $\dd$-half-disc $\dint{c}$ is not shortcut.  So $\gamma\hi\di{\db}{\dw}\geq\gamma\hi\di{\db}{\dg}$, say.  Let $\dint{a}$ be the $\dd$-disc $\di{\db}{\dg\cup\dw}$.  Then $\gamma\hi\dint{a}$ equals $\gamma\hi\di{\db}{\dg}+\gamma\hi\di{\db}{\dw}\leq n$.  Proceed to the conclusion.

\textbf{Conclusion of Cases.}
In both cases we have a $\dd$-disc $\dint{a}$ satisfying 
$\gamma\hi\dint{a}\leq n$.  Consider a $\dd$-disc $\dint{b}$ that minimizes $w\hi\dint{b}$.  Since $\dint{a}$ exists, we know $w\hi\dint{b}\leq n$.  Realize $\dint{b}$ by a geometric disc $D$ in $H$.  Then $\gamma\hi D=w\hi\dint{b}$.  Since $n$ is the girth of $\gamma$, we have $\gamma\hi D\geq n$.  Therefore 
$w\hi\dint{b}=n$.
\end{proof}

Now we present our algorithm to compute the girth of an element of a free group.  The algorithm is exponential in the rank of the free group and linear in the length of the initial word in the generators and their inverses.

\textbf{The Algorithm.}  Let $F$ be a free group of rank $g$ with generating system $\dd$.  Let $w$ be an element of $F$ represented as a word in the elements of $\dpm$.  We compute the girth of the conjugacy class of $w$.

\textbf{Step 1:  Find a shortcut $\dd$-half-disc with respect to $w$, if possible.}  There are finitely many $\dd$-half-discs.  Indeed, there are $2n(2^{2n-2}-1)$ $\dd$-half-discs, up to switching black and white.  For each $\dd$-half-disc, use Lemma~\ref{intersection free group} to check whether it is a shortcut $\dd$-half-disc with respect to $w$.  If there is a shortcut $\dd$-half-disc, apply Lemma~\ref{curve shortcut}.  Replace $\dd$ with a generating system of smaller complexity, and repeat \textbf{Step~1}.  If there is no shortcut $\dd$-half-disc, then $\dd$ is locally minimal.  Go to \textbf{Step~2}.

\textbf{Step 2:  Find a $\dd$-disc $\dint{a}$ that minimizes $w\hi\dint{a}$.}
There are finitely many $\dd$-discs.  Indeed, there are $2^{2n-1}-1$ $\dd$-discs, up to switching black and white.  For each $\dd$-disc $\dint{a}$, use Lemma~\ref{intersection free group} to calculate $w\hi\dint{a}$.  Find one that minimizes $w\hi\dint{a}$.  There may be several.  By Lemma~\ref{locally-minimal}, the quantity $w\hi\dint{a}$ equals the girth of $w$.

\textbf{End Algorithm.}

Using the results of this section, we prove here that free groups of rank two or more do not contain elements of girth one.  This gives an alternative proof of Proposition~\ref{algebraic girth one}.

\begin{lemma}
\label{group girth one}
Let $F$ be a free group.  If $F$ has rank bigger than one, then $F$ does not contain an element of girth one.
\end{lemma}

\begin{proof}
Let $F$ have rank two or more.  Let $w$ be an element of $F$ with $\girth(w)\leq 1$.  Realize $F$ as the fundamental group of a handlebody $H$.  Represent $w$ as a closed curve $\gamma$ in $H$.  Then the girth of $w$ equals the girth of $\gamma$.

Let $D$ be a disc such that $\gamma\wedge D\leq 1$.  If $\gamma\wedge D=0$, then $\gamma$ has girth zero, and we are done.  So assume $\gamma\wedge D=1$.  Homotope $\gamma$ in $H$ to minimize intersection with $D$.  By Lemma~\ref{wedge}, the curve $\gamma$ intersects $D$ exactly once.  Since $\gamma$ intersects $D$ an odd number of times, the disc $D$ is non-separating.  Complete $D$ to a disc system $\dd$.  To finish the proof, it suffices to construct a disc system of smaller complexity or show that $\gamma$ has girth zero.

Homotope $\gamma$ to minimize intersection with $\dd$.  By Corollary~\ref{subcollection}, the curve $\gamma$ continues to intersect $D$ once.  Let $D_{B}$ be one of the two oriented discs in $\dpm$ that correspond to the disc $D$.  Let $\kappa$ be the arc component of $\gamma-\dd$ that leaves $D_{B}$.  Let $D_{G}$ be the disc in $\dpm$ at the other end of $\kappa$.  Let $\di{\db}{\dw}$ denote the $\dd$-half-disc defined by taking $D_{G}$ as the gray disc $\dg$ and $D_{B}$ as the collection of black discs $\db$.  Let $\dw$ denote the white discs, which make up the rest of $\dpm$.

By construction, we have $w\wedge\di{\db}{\dw}=0$ and $w\wedge\di{\dg}{\db}=1$.  If $w\wedge\di{\dg}{\dw}=0$, then $\di{\dw}{\dg\cup\db}$ defines a $\dd$-disc $\dint{b}$ with $w\wedge\dint{b}=0$.  Realizing $\dint{b}$ as a geometric disc, we see that $w$ would have girth zero.  So assume that $w\wedge\di{\dg}{\dw}>0$.  In this case, the $\dd$-half-disc $\di{\db}{\dw}$ is a shortcut $\dd$-half-disc.  By Lemma~\ref{curve shortcut}, we can reduce the complexity of $\dd$, so we are done.
\end{proof}

\section{Main Theorem}
\label{maintheorem}
In this section we prove our main result.  In the previous section we showed how to compute the girth of an element in a free group.  The Main Theorem applies to curves in a subsurface of the boundary of a handlebody $H$.

\begin{theorem}[\textbf{Main Theorem}]
\label{subsurface theorem}
Let $H$ be a handlebody.  Let $S$ be a compact surface in $\partial H$.  The minimum girth in $H$ over all essential closed curves in $S$ can be achieved by a simple closed curve in $S$.
\end{theorem} 

We begin with some definitions and lemmas.  Let $H$ be a handlebody, and let $\gamma$ be a disjoint collection of simple closed curves in $\partial H$.  Let $\dd$ be a disc system for $H$.  Define the \emph{geometric complexity} of $\dd$ with respect to $\gamma$ to be the intersection number $|\dd\cap\gamma|$.  This parallels the notion of homotopic complexity for a word $w$ in a free group that we presented in Section~\ref{computecurve}.  We will sometimes say just \emph{complexity}.  Define a disc system to be \emph{taut} if it minimizes the complexity in its isotopy class.

\begin{figure}
\epsfig{file=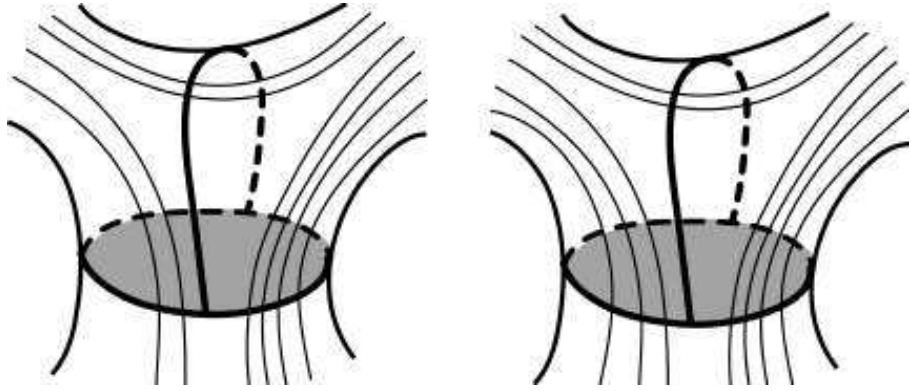, height=2in}
\caption{Both pictures illustrate a half-disc $K$ based at a shaded disc $D$.  The right one is a short-cut half-disc, and the left one is not.}
\label{halfdisc_shortcut}
\end{figure}

Let $\dd$ be a taut disc system with respect to $\gamma$.  Recall that $\dpm$ denotes the $2g$ oriented discs coming from $\dd$, where $g$ is the genus of $H$.  Define a \emph{half-disc} based at an oriented disc $D$ in $\dpm$ to be an arc $K$ in $\bd H$ with endpoints in $\bd D-\gamma$ and interior embedded essentially in $\partial H-\bd \dd$.  By essentially we mean that, inside $\bd H$, $K$ cannot be isotoped $\rel \bd K$ into $\bd \dd$.  We also require that $K$ enter and leave $D$ from the side corresponding to its orientation.  See Figure~\ref{halfdisc_shortcut} for two examples of half-discs.  Observe that, inside $H$, a half-disc $K$ can be isotoped $\rel \bd K$ into $D$.

Let $K$ be a half-disc based at a disc $D$ in $\dpm$.  The endpoints of $K$ divide $\bd D$ into two components $\bd D_B$ and $\bd D_W$.  We say that $K$ is a \emph{shortcut} half-disc if $|K\cap\gamma| < \min(|D_B\cap\gamma|,|D_W\cap\gamma|)$.  See Figure~\ref{halfdisc_shortcut} for an example.  By the following lemma, if $\dd$ has a shortcut half-disc, we can construct a new disc system with smaller complexity.  Lemma~\ref{shortcut} is the geometric analogue of Lemma~\ref{curve shortcut}.

\begin{lemma}[cf. Lemma~\ref{curve shortcut}]
\label{shortcut}
Let $H$ be a handlebody, and let $\gamma$ be a disjoint collection of simple closed curves in $\partial H$.  Let $\dd$ be a taut disc system with respect to $\gamma$.  If there is a shortcut half-disc with respect to $\gamma$, then there is a disc system with smaller complexity than that of $\dd$.
\end{lemma}

\begin{figure}
\epsfig{file=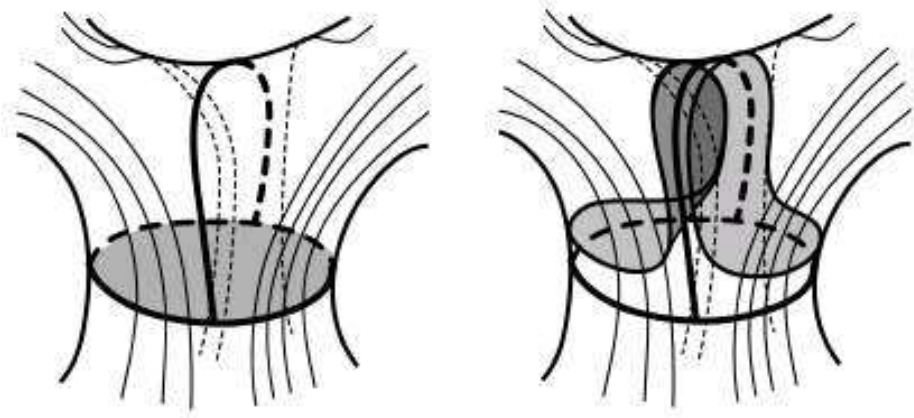, height=2.25in}
\caption{The left picture shows a compressing half-disc $K$ based at a shaded disc $D$.  The right picture shows the disc $D$ after boundary-compressing along $K$.}
\label{halfdisc_compress}
\end{figure}

\begin{proof}
Let $K$ be a shortcut half-disc based at an oriented disc 
$D$ in $\dpm$.  Let $\bd D_B$ and $\bd D_W$ be the two arc components of $\partial D-\partial K$.
Since $K$ is a shortcut half-disc with respect to $\gamma$, 
we know 
$|K\cap\gamma|<\min(|\bd D_B\cap\gamma|,|\bd D_W\cap\gamma|)$.
Boundary-compress $D$ along $K$ to create two essential 
discs in $H$ that are disjoint from $\dd$.  See Figure~\ref{halfdisc_compress}.
By the same argument used in the proof of Lemma~\ref{curve shortcut}, replacing $D$ in $\dd$ with one of these two discs, say $D^\prime$, yields a disc system.  Replace $D$ with $D^\prime$ to form a new disc system $\dd^\prime$.  Say $D^\prime$ was the disc corresponding to the arc $\bd D_B$.  Since
$$|D^\prime\cap\gamma|=|D_B\cap\gamma|+|K\cap\gamma|<
|D_B\cap\gamma|+|D_W\cap\gamma|=|D\cap\gamma|,$$ 
the disc system $\dd^\prime$ has a complexity smaller than the complexity of $\dd$.
\end{proof}

We say a disc system $\dd$ is \emph{locally minimal} with respect to $\gamma$ if there are no shortcut half-discs.  Observe that the proof of Lemma~\ref{shortcut} is constructive, so it gives an algorithm to construct a locally minimal disc system, given any collection of curves $\gamma$ in $\bd H$.

Define a \emph{compressing half-disc} to be a half-disc $K$ with $|K\cap\gamma|=0$.  See Figure~\ref{halfdisc_compress}.  We say a disc system $\dd$ is \emph{compressible} if $\dd$ has a compressing half-disc.  By the following lemma, if a curve collection $\gamma$ has a compressible locally minimal disc system, then the surface $\bd H-\gamma$ is compressible in $H$.

\begin{lemma}
\label{no waves}
Let $H$ be a handlebody, and let $\gamma$ be a disjoint collection
of simple closed curves in $\partial H$.
Let $\dd$ be a locally minimal disc system with respect to $\gamma$.  
If $\dd$ is compressible, then $H$ contains an essential disc disjoint
from $\dd$ and $\gamma$.
\end{lemma}

\begin{proof}
Let $K$ be a compressing half-disc for $\dd$ based at a disc $D$ in $\dpm$, so that $|K\cap\gamma|=0$.  Let $\bd D_B$ and $\bd D_W$ denote the two components of $\bd D-\bd K$.  Since $\dd$ is locally minimal, the half-disc $K$ is not a shortcut half-disc.  Thus $|K\cap\gamma|\geq\min(|\bd D_B\cap\gamma|,|\bd D_W\cap\gamma|)$.  In other words, at least one of the two arcs, say $\bd D_B$, satisfies $|\bd D_B\cap\gamma|=0$.  Boundary-compress $D$ along $K$ to create two essential discs in $H$ that are disjoint from $\dd$.  The disc corresponding to $\bd D_B$ can be made disjoint from $\gamma$ since neither $K$ nor $\bd D_B$ intersected $\gamma$.
\end{proof}

\begin{figure}
\epsfig{file=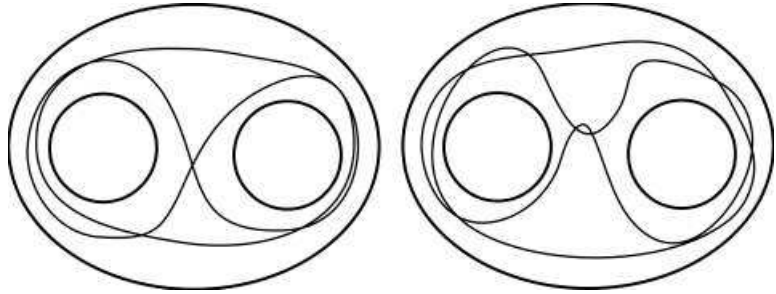, height=1.5in}
\caption{The curve in the left picture fills the surface.  The curve on the right does not fill.  It can be homotoped into an annulus.}
\label{fillingcurve}
\end{figure}

Let $S$ be a compact surface.  Let $\gamma$ be an essential closed curve in $S$.  We say that $\gamma$ \emph{fills $S$} if $\gamma$ cannot be homotoped into a subsurface $S^\prime$ of $S$ with $\ec{S^\prime}<\ec{S}$.  See Figure~\ref{fillingcurve}.  Let $S^\prime$ be a subsurface of $S$.  We say that $S^\prime$ is \emph{incompressible} in $S$ if $\pi_1(S^\prime)$ injects into $\pi_1(S)$.  By the following lemma, any curve in a surface $S$ can be homotoped to fill some incompressible subsurface in $S$.

\begin{lemma}
\label{full}
Let $S$ be a compact surface with boundary.  Let $\gamma$ be an 
essential closed curve in $S$.  Then there is an incompressible subsurface $S^\prime$ of $S$ and a curve $\gamma^\prime$ such that $\gamma^{\prime}$ both fills $S^{\prime}$ and is homotopic to $\gamma$ in $S$.
\end{lemma}

\begin{proof}
The proof is by downward induction on $\ec{S}$.  If $\ec{S}=0$, then $\gamma$ fills $S$.  Assume that $\ec{S}>0$.  If $\gamma$ does not fill $S$,
then homotope $\gamma$ into a compact subsurface $S^\prime$ of $S$ satisfying $\ec{S^\prime}<\ec{S}$.  If $S^\prime$ is compressible, then compress $S^\prime$.  The inequality $\ec{S^\prime}<\ec{S}$ still holds.  Continue compressing $S^\prime$ until it  is incompressible in $S$.  Rename $S^\prime$ as $S$, and repeat the induction.
\end{proof}

\begin{figure}
\epsfig{file=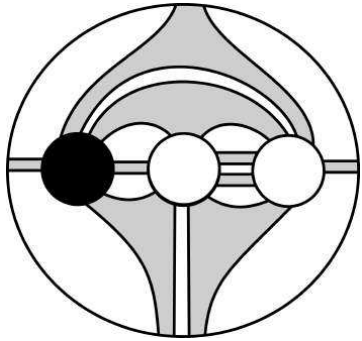, height=1.75in}
\caption{Here is a map where the handlebody has genus two and the $\dd$-interface has three white countries and one black country.
}
\label{map}
\end{figure}

Let $H$ be a handlebody of genus $g$, and let $S$ be a compact subsurface of  $\partial H$.  Let $\dd$ be a taut disc system with respect to $\partial S$.  We associate the following diagram to the triple $(H,S,\dd)$.  Take the complement in $H$ of an open regular neighborhood of $\dd$.  This defines a closed ball in $H$.  The boundary of this ball is a sphere $\m$.  The sphere $\m$ contains $2g$ distinguished discs that correspond to the discs of $\dpm$.  We call these discs the \emph{countries} of $\m$.  The surface $S$ intersected with $\m$ comprises several polygonal regions.  We call these regions the \emph{roads}.  The diagram of $\m$, together with its countries and roads, is the \emph{map} associated to $(H,S,\dd)$.  See Figure~\ref{map} for a picture.

Maps can be illustrated in the plane as follows.  Let $\m$ be a map.  View $\m$ as $\R^2\cup\infty$, the plane with a point at infinity.  Choose a country of $\m$.  Put that country around the point at infinity.  Now draw the diagram $\m$ as a disc with the other $2g-1$ countries inside.  

A $\dd$-interface adds color to a map.  Let $\m$ be a map, and let $\dint{a}$ be a $\dd$-interface.  The $\dd$-interface  $\dint{a}$ partitions $\dpm$ into black, white, and gray discs.  Color the corresponding countries of $\m$ accordingly.  Again, see Figure~\ref{map}.

\begin{figure}
\epsfig{file=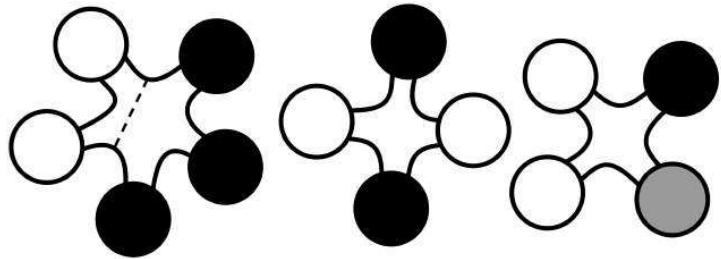, height=1.5in}
\caption{Here are three illustrations of roads and the countries they border.  The road on the left segregates.  The dotted line shown is a line of segregation.  The roads on the right do not segregate.
}
\label{segregation}
\end{figure}

Let $R$ be a road of the map $\m$.  We say the road $R$ \emph{segregates} 
if (1) it borders a black country and a white country but no gray country, and (2) a single embedded arc can separate its black borders from its white borders.  We call such an arc a \emph{line of segregation}.  See Figure~\ref{segregation}.

The following lemma relates nearly all of the concepts just introduced.

\begin{lemma}
\label{segregated regions}
Let $H$ be a handlebody, and let $S$ be a compact subsurface of $\partial H$.
Let $\dd$ be a taut, incompressible disc system with respect to $\bd S$.  Let $\dint{a}$ be a $\dd$-interface.  Let $\m$ be the associated map.  Let $\gamma$ be a closed curve in $S$, and suppose that $\gamma\wedge\dint{a}=n$.  If $\gamma$ fills $S$, then $\m$ has at most $n$ segregating roads.
\end{lemma}

\begin{proof}
Lift $\dint{a}$ to the universal cover $\hh$, and call the lift $\alpha$.  Exactly $n$ lifts of $\gamma$ cross $\alpha$.  Let $\db$ and $\dw$ denote the black and white discs of $\dpm$ associated with $\alpha$.

Since $\dd$ is incompressible with respect to $\bd S$, there are no compressing half-discs in $\bd H$.  Thus every curve in $S$ is isotopic to a geodesic with respect to $\dd$.  Isotope $\gamma$ to be a geodesic with respect to $\dd$.

View $\gamma$ in $\m$, so that $\gamma$ comprises a collection of arcs inside $\m$'s roads.  By Lemma~\ref{intersection handlebody}, every arc of $\gamma$ that connects a white country to a black country corresponds to a lift of $\gamma$ that crosses $\alpha$.

We proceed by contradiction.  Assume that $\m$ has more than $n$ segregated roads.  Then, in at least one of these roads, no arc of $\gamma$ passes from a white country to a black country.  Let $R$ be such a road.  Choose a line $\ell$ of segregation for $R$.  Isotope $\gamma$ so that it does not intersect $\ell$.  Then $S-\ell$ is a surface containing $\gamma$ and satisfying $\ec{S-\ell}<\ec{S}$.  This contradicts that $\gamma$ is full.
\end{proof}

Our final lemma is a topological result about arcs in surfaces.

\begin{lemma}
\label{arcs in S}
Let $S$ be a compact surface with boundary.  Let $n$ be a natural number such that $n\geq0$.  Let $K$ be a collection of $\ec{S}+n+1$ arcs in $S$ that are disjoint, essential, and embedded $\rel \bd S$.  If $\chi(S)\leq -1$ and $n\geq 2$, then $S$ contains an essential simple closed curve that intersects $K$ no more than $n$ times.  Otherwise, $S$ contains an essential simple closed curve that intersects $K$ no more than $n+1$ times.
\end{lemma}

\begin{figure}
\epsfig{file=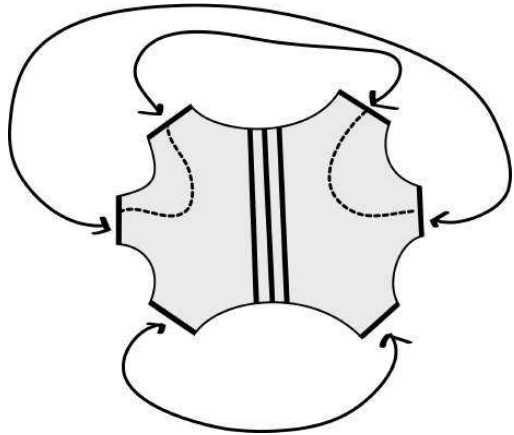, height=2.25in}
\caption{Here is an example surface $S$ for the proof of Lemma~\ref{arcs in S}.  We have $\chi(S)=-2$ and $n=3$.  Cut $S$ along $\ec{S}+1$ arcs to form the polygonal disc shown above.  Arrows indicate the gluing instructions.  The dotted line represents a curve intersecting $K$ just twice.
}
\label{polygonarcs}
\end{figure}

\begin{proof}
It is obvious when $\chi(S)=0$ that $S$ contains an essential simple closed curve that intersects $K$ exactly $n+1$ times.  So assume $\chi(S)\leq -1$.

Assume the complement of $K$ has a component that is not simply-connected.  Then there is an essential simple closed curve in $S$ that intersects $K$ zero times.  We are done in this case, so assume otherwise.

Choose a minimal collection $K_{0}$ of arcs that cut $S$ into a disc.  Any such collection contains $\ec{S}+1$ arcs.  Cut $S$ along $K_{0}$ to form a polygonal region $P$.  The polygon $P$ contains $n$ arcs leftover in its interior.  If $n=0$, there are essential curves in $S$ that intersect $K$ one time, and we would be done.  So assume that $n\geq 1$.

Let $K_{1}$ denote the $n$ interior arcs in $P$.  Equivalently, $K_{1}$ equals $K-K_{0}$.  Each arc in $K_{0}$ corresponds to a pair of sides in $\partial P$.  Consider such a pair of sides $\beta$ in $P$.  Each arc $\kappa$ in $K_{1}$ either separates or does not separate $\beta$.  An arc $\kappa$ is said to \emph{separate} $\beta$ if and only if $\beta$ cannot be connected by a path in $P$ that is disjoint from $\kappa$.

Assume there is an arc $\kappa$ in $K_{1}$ that does not separate some pair $\beta$.  Then $\beta$ can be connected by a path in $P$ that intersects no more than $n-1$ arcs in $K_{1}$.  In other words, the surface $S$ contains an essential simple closed curve that intersects no more than $n$ arcs in $K$.  In this case, we are done.  So assume otherwise, namely that all arcs in $K_{1}$ separate all pairs $\beta$.

For every pair $\beta$ in $P$, all paths connecting $\beta$ in $P$ intersect all $n$ arcs in $K_{1}$.  This implies that the $n$ interior arcs $K_{1}$ are parallel, as shown in Figure~\ref{polygonarcs}.  Since $S$ is not annulus, we can find an essential simple closed curve in $S$ that intersects just two arcs in $K_{0}$ and avoids all $n$ arcs in $K_{1}$.  Again, see Figure~\ref{polygonarcs} for an illustration.  Since $S$ contains a curve intersecting $K$ exactly two times, the proof is complete for both $n=1$ and $n\geq 2$.
\end{proof}

We now prove the Main Theorem.  We state and prove 
Lemmas~\ref{main lemma easy} and~\ref{main lemma corollary} 
in the next section, Section~\ref{mainlemma}.

\begin{proof}[Proof of Theorem~\ref{subsurface theorem}]
Assume without loss of generality that $S$ is incompressible in $\bd H$.  Then no component of $\bd S$ bounds a disc in $\bd H$.

Let $\gamma$ be a curve in $S$ that minimizes the girth in $H$ over all essential closed curves in $S$.  Let $n$ be the girth of $\gamma$.  Apply Lemma~\ref{full}, and homotope $\gamma$ into an incompressible subsurface $S^\prime$ of $S$ so that $\gamma$ fills $S^\prime$.

It suffices to prove the theorem for $S^\prime$.  Give $S^\prime$ the name $S$ for convenience.  If $S$ is an annulus, then we are done by Proposition~\ref{power}: a core curve for the annulus minimizes the girth.  So assume $S$ is not an annulus.  In particular, we can also assume $\genus(H)\geq2$.

The curve $\gamma$ fills the surface $S$.  We will construct an essential disc $D$ in $H$ and an essential simple closed curve $\gamma^\prime$ in $S$ so that $\gamma^\prime$ intersects $D$ in $n$ points.

Using Lemma~\ref{shortcut}, find a locally minimal disc system $\dd$ with respect to $\bd S$.  Then by Lemma~\ref{no waves}, either $\dd$ is incompressible with respect to $\bd S$, or there is a disc $D$ in $H-\dd$ disjoint from $\partial S$.  In the latter case, $\bd S$ contains simple closed curves of girth zero, and we are done.  So assume that $\dd$ is incompressible.

Let $B$ be a girth-realizing disc for $\gamma$, so that $\gamma\wedge B=n$.  Apply Lemma~\ref{subordination-lemma} to $B$.  There are two cases.

\textbf{Case 1:  There is a $\dd$-disc $\dint{b}$ subordinate to $B$.}
Then $\gamma\wedge\dint{b}\leq n$.  By Lemma~\ref{segregated regions}, the $\dd$-disc $\dint{b}$ has no more than $n$ segregated roads.  Apply Lemma~\ref{main lemma easy} to obtain a disc $D$ in $H-\dd$ that intersects $S$ in no more than $\ec{S}+n+1$ arcs.  Now proceed to the conclusion following Case 2.

\textbf{Case 2:  There are two $\dd$-half-discs independently subordinate
to $A$.}  Let $\dint{c}$ be the $\dd$-half-disc that $\gamma$ homotopically 
intersects less.  Let $m$ equal $\gamma\wedge\dint{c}$.  Then $m\leq\lfloor\frac{n}{2}\rfloor$.  By Lemma~\ref{segregated regions}, the $\dd$-half-disc $\dint{c}$ has no more than $m$ segregated roads.
Apply Lemma~\ref{main lemma corollary} to obtain a disc $D$ in $H-\dd$
that intersects $S$ in no more than $\ec{S}+2m+1\leq\ec{S}+n+1$ arcs.
Proceed now to the conclusion below.

\textbf{Conclusion of Cases.}
So far we have constructed a disc $D$ in $H$ that intersects $S$ in no more than $\ec{S}+n+1$ arcs.  Since $S$ is not an annulus, we have $\ec{S}\leq-1$.  Apply Lemma~\ref{arcs in S} to $S$.  There are two possible conclusions.

If $n=0$, Lemma~\ref{arcs in S} gives us an essential simple closed curve $\gamma^{\prime}$ in $S$ that intersects $D$ no more than once.  If necessary, apply Lemma~\ref{geometric girth one} to construct a disc that intersects $\gamma^\prime$ zero times.  This completes the case of $n=0$.

Assume next that $n\geq1$.  Since $n$ is the girth of a curve, we have by Lemma~\ref{group girth one} that $n\neq1$.  Hence $n\geq2$.  In this case, Lemma~\ref{arcs in S} gives us an essential simple closed curve $\gamma^{\prime}$ in $S$ that intersects $D$ no more than $n$ times.  Thus we are done.
\end{proof}

\section{Main Lemma}
\label{mainlemma}
In the previous section we proved our main theorem, Theorem~\ref{subsurface theorem}, except for two key lemmas.  We state and prove those lemmas in this section.

Lemmas~\ref{main lemma easy} and~\ref{main lemma corollary} apply to Case~1 and Case~2 of our proof of Theorem~\ref{subsurface theorem}, respectively.  Lemma~\ref{main lemma easy} applies to the $\dd$-disc case, and Lemma~\ref{main lemma corollary} applies to the $\dd$-half-disc case.  Lemma~\ref{main lemma easy} is easy to prove, and Lemma~\ref{main lemma corollary} is harder.  We present Lemma~\ref{main lemma corollary} as a corollary of Lemma~\ref{main lemma}, the Main Lemma.  We introduce several definitions before proceeding to the proofs.

Let $H$ be a handlebody.  Let $S$ be a compact surface in $\partial H$.  Let $\dd$ be a disc system in $H$.  Assume that $\dd$ is taut with respect to $\partial S$.  We defined taut in Section~\ref{maintheorem}.  Also in that section, we defined the map $\m$ associated to the triple $(H, S, \dd)$.  The map $\m$ is a spherical diagram showing the roads $R$ and the collection $\dpm$ of countries.  A $\dd$-interface $\dint{a}$ for $\dpm$ colors each country in $\m$ either black, white, or gray.

In Section~\ref{maintheorem}, we explained how to view $\m$ as a disc by choosing one country to be the outside.  This lets us refer to the inside and outside of simple closed curves in $\m$.  In this section, we will always view maps $\m$ in this way.  We denote the boundary of the outside country $\bd\m$.  When $\dint{a}$ is a $\dd$-half-disc, we choose the gray disc to be the outside country.  Otherwise, our convention is to pick a black disc as the outside.

Every road $R$ in $\m$ has $p$ \emph{prongs}, where $2p$ is the number of sides of $R$.  For example, the roads in Figure~\ref{segregation} have five, four, and four prongs, respectively.  Each road in $\m$ has at least two prongs.  Define the Euler characteristic $\chi(R)$ of a road $R$ to be $\frac{2-p}{2}$.  For example, a rectangular road has zero Euler characteristic.  All other roads in $\m$ have negative Euler characteristic.  Define the Euler characteristic of a collection of roads to be the sum of the Euler characteristics of the roads.

\begin{figure}
\epsfig{file=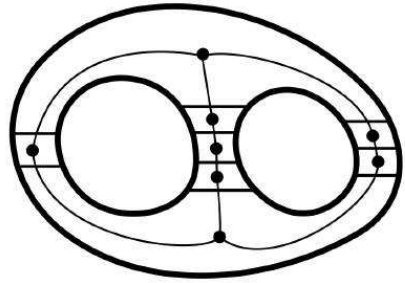, height=1.5in}
\caption{
Each 3-pronged road in this picture contributes $-1/2$ to the Euler characteristic of $S$, for a total Euler characteristic of $-1$.  The dual graph also has Euler characteristic $-1$.  Each 3-pronged vertex contributes $-1/2$ to the Euler characteristic of the graph.
}
\label{dual graph}
\end{figure}

Observe that the Euler characteristic of the collection of all roads in $\m$ equals $\chi(S)$.  This can be seen by considering the graph dual to $S$.  View the roads inside $S$ rather than inside $\m$.  Put a vertex inside each road of $S$.  For each place where two roads border, connect the corresponding vertices with an edge.  See Figure~\ref{dual graph} for an example.  The dual graph has the same Euler characteristic as the surface $S$ because the dual graph is a deformation retract of $S$.  The contribution near a vertex to the Euler characteristic of the graph equals the contribution of the corresponding road to $\chi(S)$.

Let $R$ be a road of $\m$.  Define $d(R)$ to be $\min(1,\ec{R})$.  Then $d(R)$ is $0$, $1/2$, or $1$ depending on whether $R$ has two, three, or more than three prongs.  Let $X$ be a subset of $\m$.  Recall the definition of a non-segregating road in Section~\ref{maintheorem}.  Define $d(X)$ to be the sum of $d(R)$ over all \emph{non-segregating} roads $R$ that intersect $X$.  Loosely speaking, the function $d$ counts the non-segregating roads of negative curvature that intersect the subset $X$ of $\m$.  A non-segregating road $R$ with $d(R)$ equal $1/2$ must be a 3-pronged road touching a gray $\bd\m$.  Observe that any subset $X$ of $\m$ satisfies $d(X)\leq\ec{S}$.

\begin{figure}
\epsfig{file=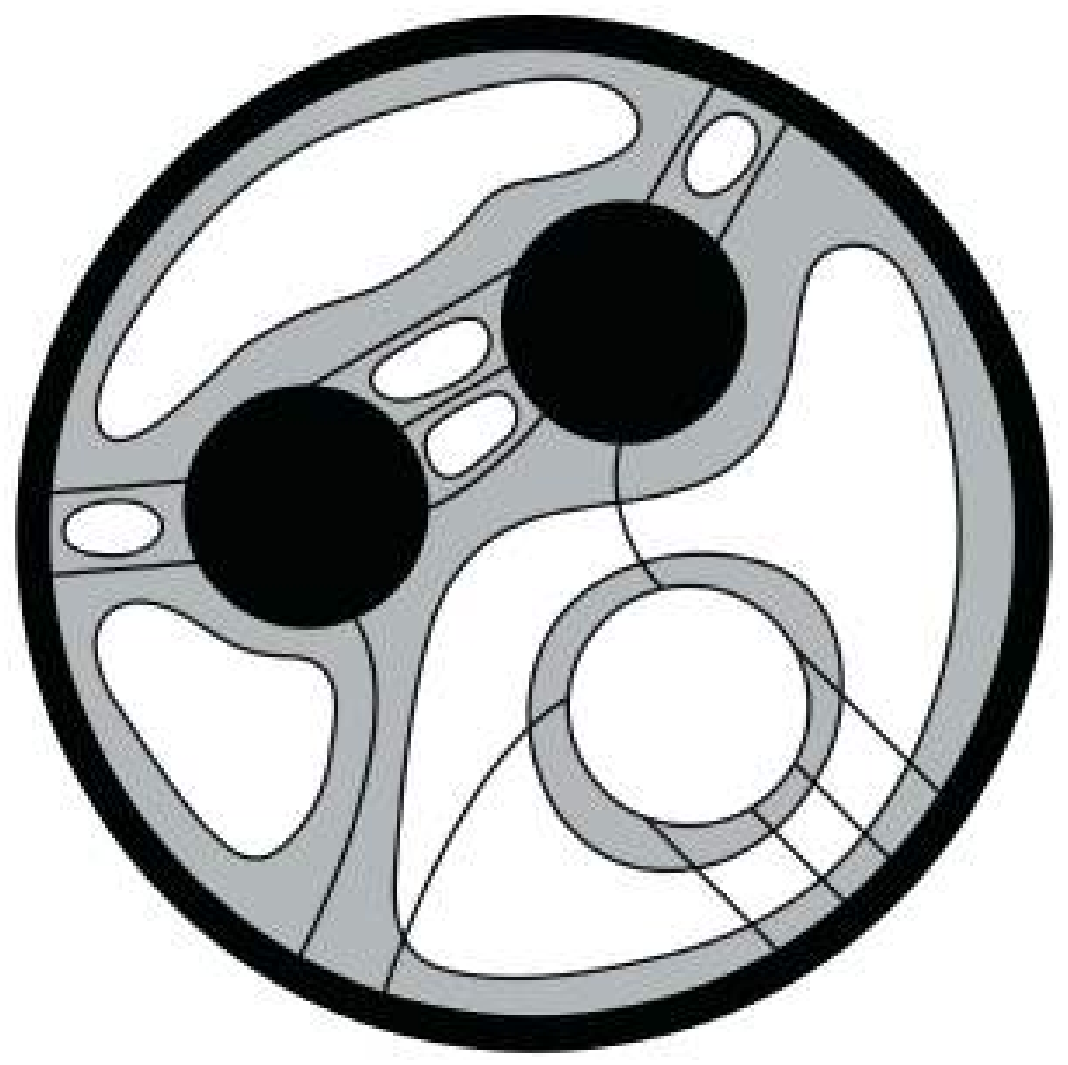, height=2.25in}%
\hspace{.5in}%
\epsfig{file=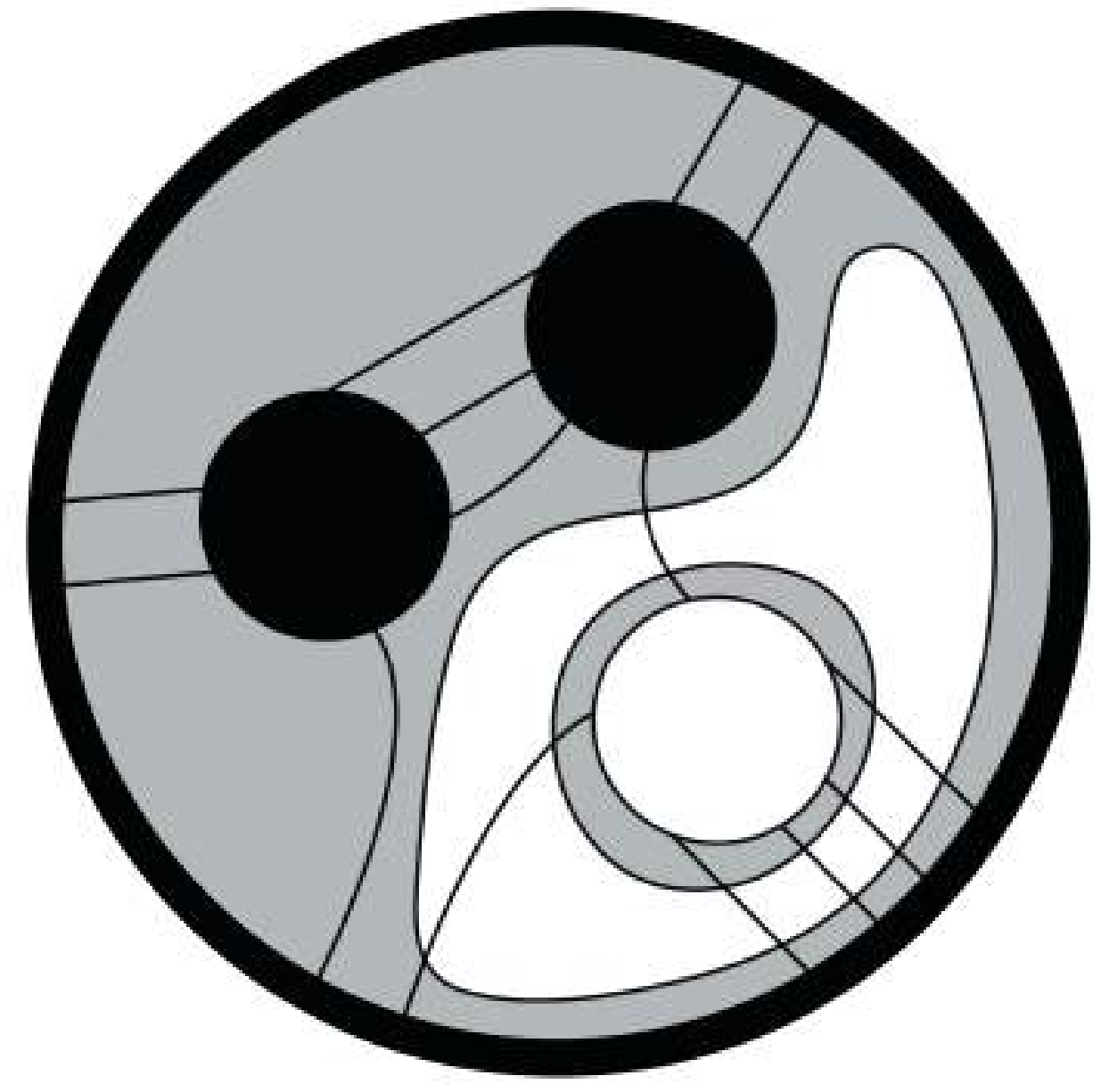, height=2.25in}
\caption{
We show two pictures of the same map $\m$ with a $\dd$-disc.  The outer country is black.  We show the collection $U$ of continents in gray.  The left picture shows $U$ before compressing.  The right shows $U$ after compressing.  There is a black continent and a white continent.
}
\label{map1}
\end{figure}

\begin{figure}
\epsfig{file=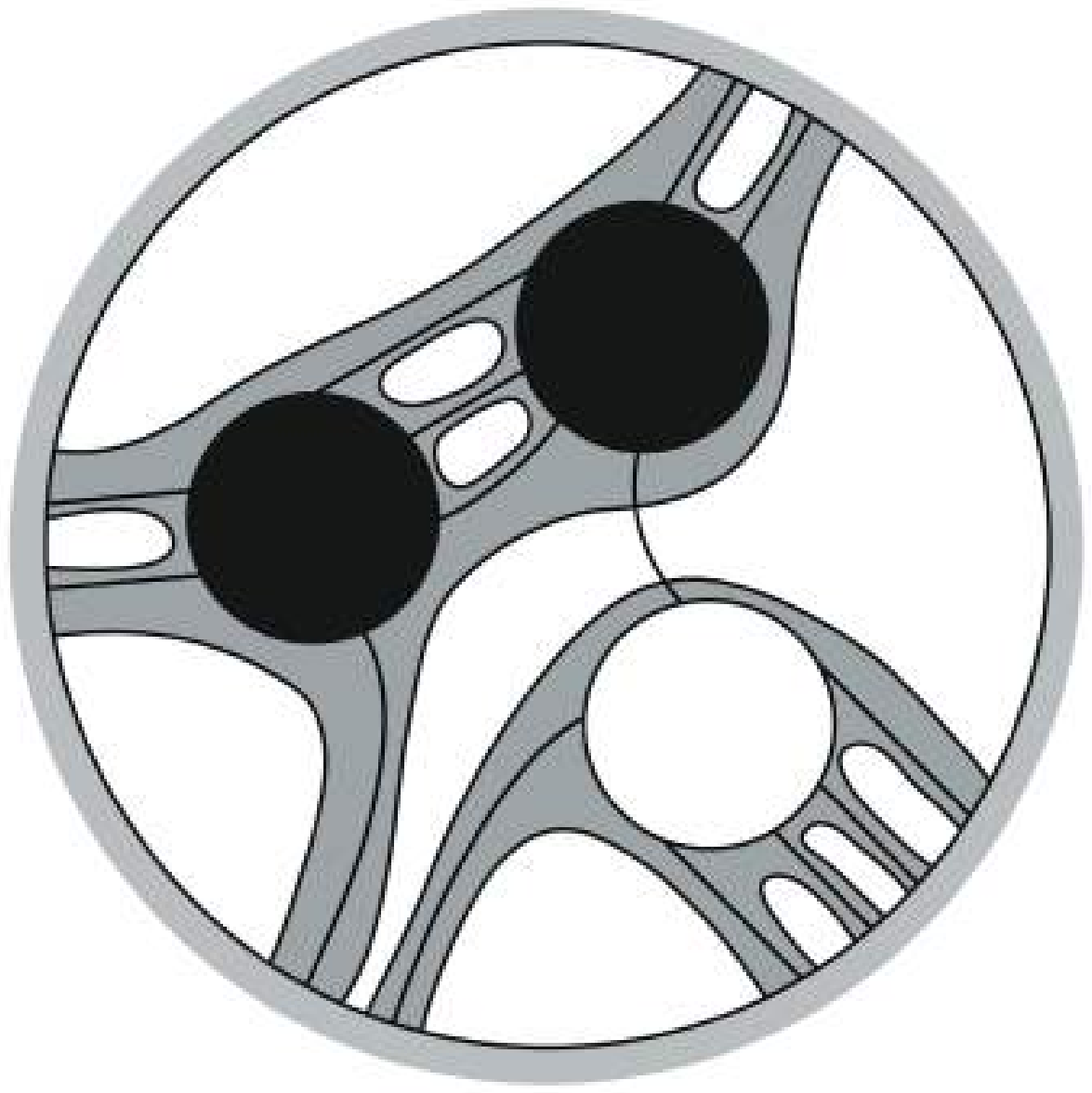, height=2.25in}
\hspace{.5in}
\epsfig{file=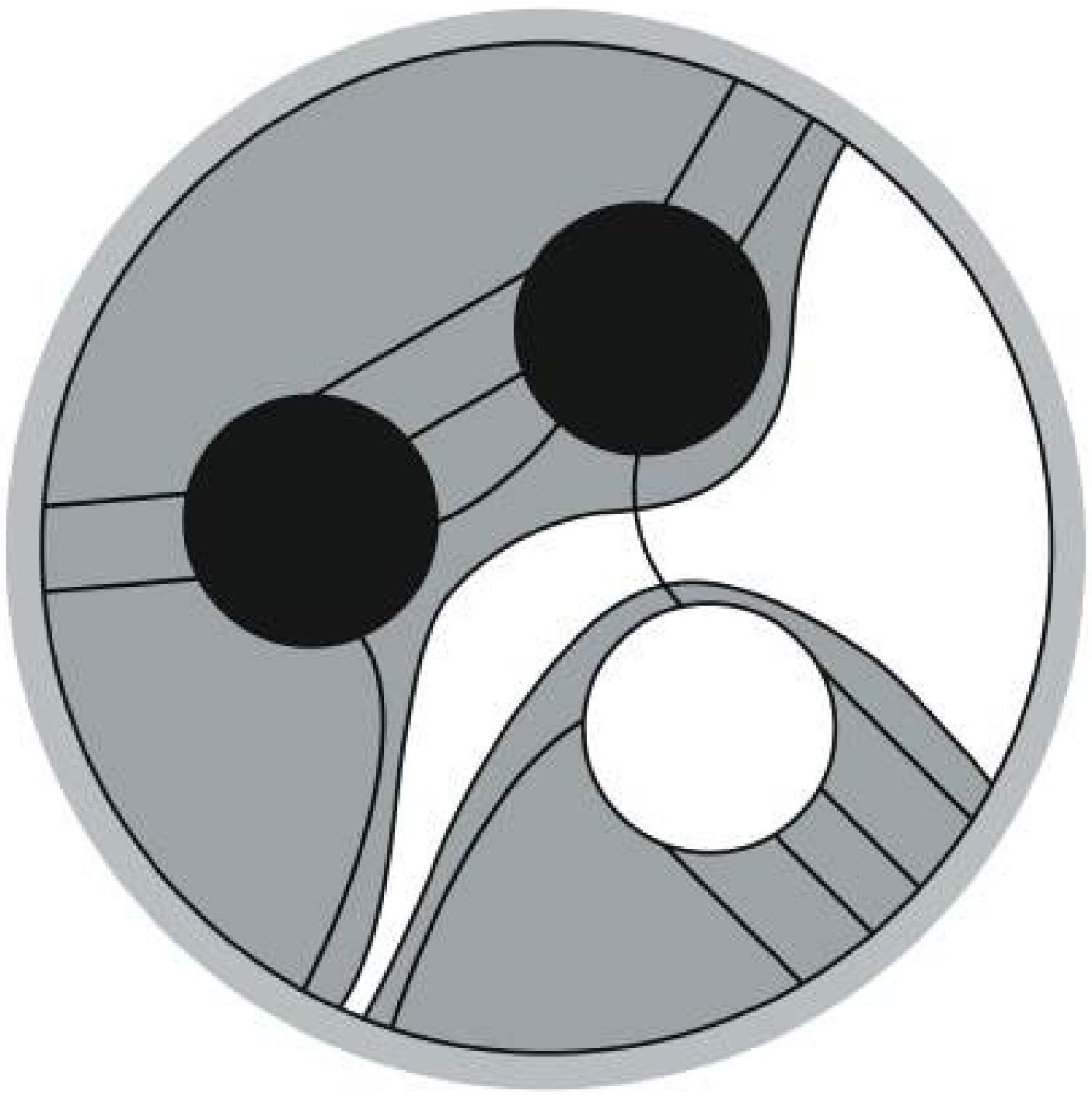, height=2.25in}
\caption{
We show a map $\m$ with a $\dd$-half-disc.  The boundary $\bd\m$ corresponds to a gray country.  The left shows the collection $U$ of continents before compressing, and the right after compressing.  Again, there is a black continent and a white continent.}
\label{graymap}
\end{figure}

We now define the continents of $\m$.  Consider the union in $\m$ of (1) all the countries in $\m$, (2) the boundary circle $\partial \m$ if it corresponds to a black country, and (3) all the arcs of $\partial S$ that do not connect a black country to a white country.  Let $U$ be a regular neighborhood of this subset in $\m$.  See the left sides of Figures~\ref{map1} and~\ref{graymap} for examples.  Each component of $\partial U$ is either a circle or an arc.  Arc components exist only when $\dint{a}$ is a $\dd$-half-disc (see Figure~\ref{graymap}).  In this case, all arcs begin and end at the boundary circle $\partial \m$.

Compress the inessential boundary components of $U$
as follows.
Fill in each circle of $\partial U$ that bounds a disc in $\m-\dd$.
Boundary-compress every arc of $\partial U$ that is inessential
in $\m-\dd$.
Each component of $\bd U$ is now either an essential circle in 
$\m-\dd$ or an essential arc in $\m-\dd$ based at $\partial \m$.  
See the right sides of Figures~\ref{map1} and~\ref{graymap} for
examples.

A \emph{continent} of $\m$ is a connected component of the
neighborhood $U$ after compressing.  Each continent 
contains at least one country.
Each continent contains only like-colored countries.
Thus we can speak of black and white continents.

\begin{lemma}
\label{main lemma fact}
Let $U$ be the collection of continents for a $\dd$-interface
$\dint{a}$, as defined above.  Let $m$ be the number of 
segregating roads.  
Let $K$ be a component of $\partial U$ that intersects each road 
in at most one arc.  Then $|K\cap \partial S|\leq 2d(K)+2m$.
\end{lemma}

\begin{proof}
The curve $K$ intersects the boundary of a non-segregating road with 
two prongs in zero points, with three prongs in at most one point, 
and with four or more prongs in at most two points.
The curve $K$ also intersects every segregated road in at most 
two points.  The result follows.
\end{proof}

\begin{figure}
\epsfig{file=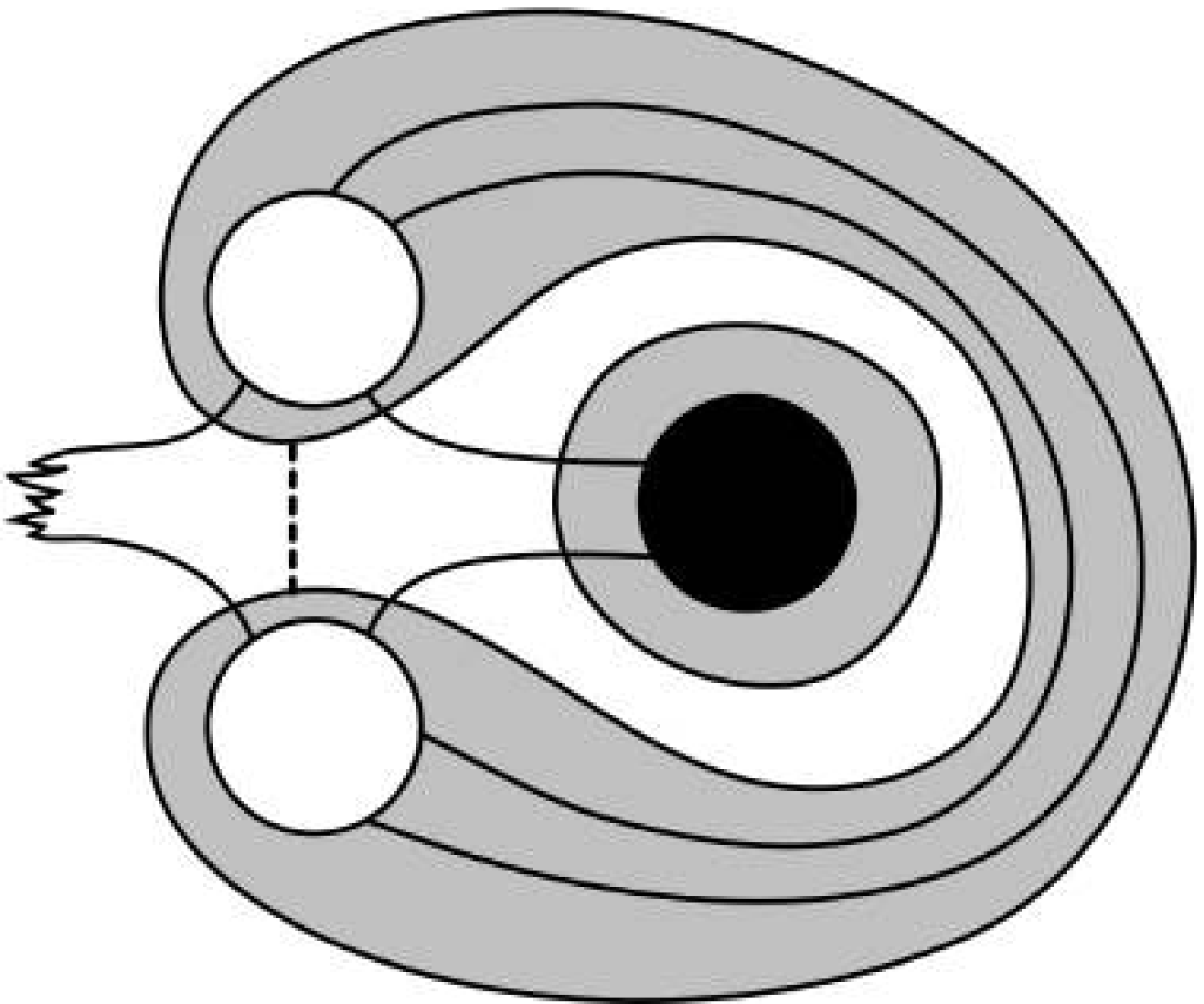, height=2.25in}%
\caption{
These two components $K_1$ and $K_2$ of $\bd U$ satisfy 
$K_1\subseteq K_2$.  The curve $K_1$ lies inside the union
of $K_2$ and the dotted arc.
}
\label{partialorder}
\end{figure}

We now define a partial order $\subseteq$ on the components of 
$\partial U$.  Let $K_1$ and $K_2$ be distinct components of 
$\partial U$.  We say $K_1\subseteq K_2$ if there is a
simple arc $k$ in some road so that $K_1$ is inside $K_2\cup k$
(see Figure~\ref{partialorder}).
We also define $K_1\subseteq K_1$ for all components $K_1$
of $\bd U$.  This makes $\subseteq$ a partial order on the components
of $\partial U$.

Observe that if $K_1\subseteq K_2$ for distinct $K_1$ and $K_2$, then $K_1$ is a circle rather than an arc.  Let $K_2$ be a component of $\partial U$.  Observe that there is a distinct component $K_1$ of $\partial U$ such that $K_1\subseteq K_2$ if and only if $K_2$ intersects a road in more than one arc.  In particular, the boundary $\partial U$ contains a circle.

We now state Lemma~\ref{main lemma easy}.
The proof we give specializes our proof of the Main Lemma~\ref{main lemma}.

\begin{lemma}[cf. Lemma~\ref{main lemma corollary}]
\label{main lemma easy}
Let $H$ be a handlebody, and let $S$ be a compact surface in 
$\partial H$.  Let $\dd$ be a taut disc system with respect 
to $\partial S$.  Let $\dint{b}$ be a $\dd$-disc with 
$n$ segregated roads.
Then there is an essential disc in $H-\dd$ that intersects $S$ in no more
than $\ec{S}+n$ arcs.
\end{lemma}

\begin{proof}
Consider the map $\m$ corresponding to $(H,S,\dd)$, with 
$\dd$-disc $\dint{b}$.
Let $U$ denote the union of continents in $\m$, as defined 
above.  Since no countries are gray, all components of 
$\partial U$ are circles.  Let $C$ be a circle that is minimal 
with respect to the partial ordering $\subseteq$.  Then $C$ 
intersects each road in at most one arc.
Apply Lemma~\ref{main lemma fact} to $C$.  Then 
$|C\cap\partial S|\leq 2d(C)+2n\leq 2\ec{S}+2n$.  Thus
$C$ intersects $S$ in at most $\ec{S}+n$ arcs.  Since $C$
is disjoint from $\dd$ in $H$, it bounds a disc in $H$.
\end{proof}

We now state and prove the Main Lemma.

\begin{lemma}[\textbf{Main Lemma}]
\label{main lemma}
Let $H$ be a handlebody, and let $S$ be a compact surface in $\partial H$.  Let $\dd$ be a taut disc system with respect to $\partial S$.  Let $\dint{c}$ be a $\dd$-half-disc with $m$ segregating roads. Then in $H-\dd$ there is either (1) an essential circle $C$ such that $|C\cap\partial S|\leq2(\ec{S}+2m+1)$, or (2) a half-disc $K$ such that $|K\cap\partial S|\leq\ec{S}+2m+1$.
\end{lemma}

\begin{proof}
Let $\m$ be the associated map.  The $\dd$-interface $\dint{c}$ gives a color to each disc in $\dpm$.  Since $\dint{c}$ is a $\dd$-half-disc, exactly one disc in $\dpm$ is gray.  The discs in $\dpm$ correspond to the countries of $\m$.  Following our convention, the disc corresponding to $\bd \m$ is the gray disc.

Let $U$ denote the union of continents in $\m$.  Assume temporarily that $\partial U$ contains a circle.  Let $C$ be a circle in $\bd U$ that is minimal with respect to the partial ordering $\subseteq$.  Then $C$ intersects each road in at most one arc.  Recall this by examining Figure~\ref{partialorder}.  Applying Lemma~\ref{main lemma fact}, we have $|C\cap\partial S|\leq 2d(C)+2m\leq 2\ec{S}+2m\leq 2(\ec{S}+2m+1)$.  We are done in this case, so we can assume otherwise.

\textbf{Assume that every component of $\bd U$ is an arc.}  Then every component $K$ of $\bd U$ intersects each road in at most one arc.  We will construct an arc $K$ of $\bd U$ that satisfies the conclusion of the lemma.  

Give the complement of $U$ in the disc $\m$ the structure of a graph, as follows.  The disc $\m$ contains a collection of arcs coming from $\bd S$.  These arcs bound the roads in $\m-\dpm$.  Let each component of the complement of these arcs in $\m-U$ be a vertex.  Note that some of these components correspond to roads, and others do not.  Let each arc in $\m-U$ be an edge.  The incidence relation gives these vertices and edges the structure of a finite graph.  Observe that each road in $\m$ corresponds to only one vertex in the graph.  Note also that this graph may be disconnected.

Let $\Gamma$ be a connected component of this graph.  Since $\bd U$ contained no circles, the graph $\Gamma$ is a tree.  Observe that every second vertex of $\Gamma$ corresponds to a road of $\m$.  More precisely, let $v$ be a vertex of $\Gamma$ that corresponds to a road.  The vertices of $\Gamma$ that correspond to roads are exactly the vertices that are an even number of edges away from $v$.


Define the \emph{extremal} vertices of a tree to be the vertices of the tree that touch only one edge.  Define the \emph{interior} vertices of a tree to be the non-extremal vertices.  Observe that if an extremal vertex of the tree $\Gamma$ corresponds to a road, then that vertex corresponds to a road touching $\bd\m$ with three or more prongs.  Note also that roads touching $\bd\m$ are non-segregating. 

Assume temporarily that no interior vertex of $\Gamma$ corresponds to a non-segregating road.  Then every interior vertex that corresponds to a road corresponds to a segregating road.  Choose a geodesic path $P$ in $\Gamma$ that connects two extremal vertices.  Realize $P$ as a half-disc $K$ in $\m$ based at the gray disc corresponding to $\bd\m$.  Observe that $|K\cap S|$ equals the number of edges of $P$.  There are two cases.

\textbf{Case 1: Neither endpoint of $P$ corresponds to a road.}  Then the number of edges of $P$ is less than or equal to twice the number of segregating roads in $\m$.  Hence $|K\cap S|\leq 2m<\ec{S}+2m+1$, and we are done with this case. 

\textbf{Case 2. At least one endpoint of $P$ corresponds to a road.}  Then at least one road in $\m$ has three or more prongs, so $\ec{S}>0$, or equivalently $\ec{S}\geq 1$.  The number of edges of $P$ is less than or equal to $2m+2$, since each endpoint of $P$ can contribute an additional edge.   Hence $|K\cap S|\leq 2m+2\leq\ec{S}+2m+1$.  We are done with the second case and with our temporary assumption. 

\textbf{Assume that at least one interior vertex of $\Gamma$ corresponds to a non-segregating road.}  Such an interior vertex must correspond to a road with four or more prongs, so $\ec{S}\geq 1$.

\begin{figure}
\epsfig{file=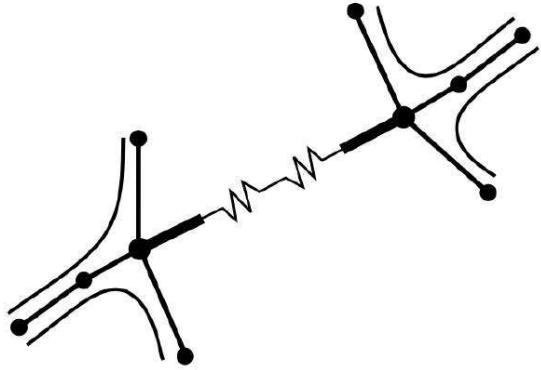, height=2in}
\caption{
We illustrate how four paths can be chosen in the graph $\Gamma$ to satisfy the desired properties.  The subgraph $\Gamma_{0}$ is shown in bold, though much of the graph $\Gamma$ is not shown.  Each of the four paths intersects exactly one extremal vertex of $\Gamma_{0}$.  No more than two of the four paths intersect any one extremal vertex of $\Gamma$.
}
\label{four paths}
\end{figure}

Let $\Gamma_{0}$ denote the convex hull in $\Gamma$ of the interior vertices that correspond to non-segregating roads.  By assumption this is nonempty.  Choose four geodesic paths $P_{1}$, $P_{2}$, $P_{3}$, and $P_{4}$ in $\Gamma$ such that (1) no more than two of the $P_{i}$ touch any extremal vertex of $\Gamma$, and (2) each $P_{i}$ intersects $\Gamma_{0}$ in at most one vertex.  This can be done easily.  The process is illustrated in Figure~\ref{four paths}.

The argument now divides into cases according to how many extremal vertices of $\Gamma$ correspond to roads.  The argument is similar to above.  Note that if distinct extremal vertices of $\Gamma$ correspond to roads, they correspond to distinct roads.

Represent each geodesic path $P_{i}$ in $\Gamma$ by a half-disc $K_{i}$ in $\m$ based at the gray disc $\bd\m$.  Observe that $|K_{i}\cap S|$ equals the number of edges of $P_{i}$.  At most four edges of each $P_{i}$ do not correspond to the intersection of $K_{i}$ with the boundary of a segregating road: possibly two edges for an extremal vertex of $\Gamma_{0}$ and possibly one edge for each extremal vertex of $P_{i}$.  There are three cases.

\textbf{Case 1. Zero or one extremal vertices of $\Gamma$ correspond to roads.}
Choose one of the four paths, call it $P_{1}$, that does not intersect an extremal vertex of $\Gamma$ that corresponds to a road.  Then at most two edges of $P_{1}$ correspond to non-segregating roads.  The rest correspond to segregating roads.  Thus we have $|K_{1}\cap S|\leq 2m+2 \leq \ec{S}+2m+1$.  The latter inequality holds since $\ec{S}\geq 1$.

\textbf{Case 2. Two or three extremal vertices of $\Gamma$ correspond to roads.}
Then at least one of the four paths, call it $P_{1}$, intersects at most one extremal vertex of $\Gamma$ that corresponds to a road.  Such a road must be non-segregating and contributes one additional edge.  Moreover, each extremal vertex contributes at least $1/2$ to $\ec{S}$, so that $\ec{S}\geq 2$.  Thus $|K_{1}\cap S|\leq 2m+3 \leq \ec{S}+2m+1$.

\textbf{Case 3. Four or more extremal vertices of $\Gamma$ correspond to roads.}
Then each of the four paths intersect at most two extremal vertices of $\Gamma$ corresponding to a road.  This contributes at most two additional edges.  Moreover, each extremal vertex contributes at least $1/2$ to $\ec{S}$, so that $\ec{S}\geq 3$.  Taking $P_{1}$ for example, we have $|K_{1}\cap S|\leq 2m+4 \leq \ec{S}+2m+1$.

In each case we found a half-disc $K$ such that $|K\cap S|\leq \ec{S}+2m+1$, so we are done.
\end{proof}

\begin{lemma}[corollary of Main Lemma~\ref{main lemma}; 
cf. Lemma~\ref{main lemma easy}]
\label{main lemma corollary}
Let $H$ be a handlebody, and let $S$ be a compact surface in $\partial H$.
Let $\dd$ be a locally minimal disc system with respect to $\partial S$.
Let $\dint{c}$ be a $\dd$-half-disc with $m$ segregated roads.
Then there is an essential disc in $H-\dd$ that intersects $S$ in no more than
$\ec{S}+2m+1$ arcs.
\end{lemma}

\begin{proof}
Let $D$ be the disc in $\dpm$ at which the $\dd$-half-disc $\dint{c}$ is based.  Apply Lemma~\ref{main lemma} to $\dint{c}$.  There are two cases.

\textbf{Case 1.}  There is an essential circle $C$ in $\m-\dpm$ such 
that $|C\cap\partial S|\leq2\ec{S}|+4m+2$.  Then $C$ intersects
$S$ in no more than $\ec{S}+2m+1$ arcs.  The circle $C$ bounds
a disc in $H-\dd$, so we are done.

\textbf{Case 2.}There is a half-disc $K$ based at $D$ such that $|K\cap\partial S|\leq\ec{S}+2m+1$.  The endpoints $\partial K$ divide $\partial D$ into two components $D_B$ and $D_W$.  Since $\dd$ is locally minimal, the half-disc $K$ is not a short-cut
half-disc.  Therefore $D_B$, say, satisfies $|D_B\cap\partial S|\leq |K\cap\partial S|$.  The circle $C$ formed by $D_B$ and $K$ is essential and satisfies $|C\cap\partial S|\leq 2\ec{S}+4m+2$, so we can finish as in Case~1.
\end{proof}

\section{Computing the Girth of a Subsurface}
\label{surfacealgorithm}
\label{computesurface}
In the previous two sections we proved the Main Theorem.  We showed that every compact surface in the boundary of a handlebody contains a simple closed curve that minimizes the girth over all closed curves in the surface.  In this section we describe how to construct such a curve.  This gives us an algorithm to compute the girth of any surface.  To describe the algorithm, we collect together some of the key constructions found in earlier sections of this paper.

Let $H$ be a handlebody, and let $S$ be a surface contained in 
$\bd H$.  Let $\dd$ be a disc system for $H$.  Recall from 
Section~\ref{introduction} that a disc system for $H$ is a disjoint
collection of discs that divides $H$ into a ball.  In 
Section~\ref{maintheorem}, we described a way to simplify
a disc system with respect to $S$.  Lemma~\ref{shortcut} describes 
the procedure.  It involves boundary-compressing along certain half-discs
until a locally minimal disc system is reached.
Let \locallyminimize{$\dd$} denote a disc system obtained 
from $\dd$ using this procedure.

Let $\dint{a}$ be a $\dd$-disc or a $\dd$-half-disc.  
Recall the associated diagram $T$.  Each region in $T$ is 
either segregated or desegregated.  
Let \segregation{$\dint{a}$} denote the number of regions in $T$ that
are segregated with respect to $\dint{a}$.

Let $n$ be \segregation{$\dint{a}$}.
When $\dint{a}$ is a $\dd$-disc (resp. $\dd$-half-disc),
Lemma~\ref{main lemma easy} (resp. Lemma~\ref{main lemma corollary})
describes a procedure to find a disc intersecting $S$ in no more
than $\ec{S}+n$ arcs (resp. $\ec{S}+\max(1,2n)$ arcs).
Let \constructdisc{$\dint{a}$} denote a disc obtained from $\dint{a}$ using 
this procedure.

Let $D$ be a disc intersecting $S$ in no more than $\ec{S}+n$ arcs.
Lemma~\ref{arcs in S} describes a procedure to find a simple
closed curve in $S$ that intersects $D$ in no more than $\max(1,n-1)$ 
points, or $n$ points if $\chi(S)$ is zero.  
Let \constructcurve{$D$} denote a simple closed curve obtained from $D$
using this procedure.

\subsection{The Algorithm.}
\label{surface algorithm}
Let $H$ be a handlebody, and let $S$ be a surface in $\partial H$.  
The algorithm below can be used to find a girth-minimizing curve in $S$
and a girth-realizing disc.

The algorithm takes as input a natural number $n$, a subsurface $S$,
and a disc system $\dd$.  It either terminates with no output or 
terminates after returning a simple closed curve $\gamma$ in $S$ and 
an essential disc $D$.  If it terminates with no output, then $\girth(S)>n$.
If it terminates with output, then $\girth(S)\leq n$.  The curve and disc 
satisfy $|\gamma\cap D|\leq\max(1,n)$.
Thus $\girth(\gamma)\leq n$.  The disc $D$ is girth-realizing for 
$\gamma$ unless they intersect once.
In this case a new disc can easily be formed that is disjoint from
$\gamma$ (Proposition~\ref{geometric girth one}).

To compute the girth, choose an arbitrary curve $\gamma$ in $S$.  
Compute the girth of $\gamma$ using the algorithm of 
Section~\ref{curvealgorithm}, say.  This establishes an upper bound
$n$ on the girth of $S$.  Repeat the algorithm below for successively
smaller $n$ until it terminates with no output.

The algorithm below is recursive in $\ec{S}$.  
The function \texttt{STOP} terminates all computation.

\bigskip

\begin{ttfamily}
Input: natural number $n$, surface $S$, disc system $\dd$

\smallskip

begin \GIRTH{$n$,$S$,$\dd$};

\hspace{20 pt}if ( $\chi(S)=2$ ) end GIRTH($n$,$S$,$\dd$);

\hspace{20 pt}$\dd\leftarrow$\locallyminimize{$\dd$};

\hspace{20 pt}for each $\dd$-disc $\dint{a}$;

\hspace{40 pt}if \segregation{$\dint{a}$}$\leq n$;

\hspace{60 pt}$D=$\constructdisc{$\dint{a}$};

\hspace{60 pt}$\gamma=$\constructcurve{$D$};

\hspace{60 pt}print $(\gamma, D)$;

\hspace{60 pt}STOP;

\hspace{40 pt}for each line of segregation $\ell$;

\hspace{60 pt}\GIRTH{$n$,$S-\ell$,$\dd$};

\hspace{20 pt}for each $\dd$-half-disc $\kappa$;

\hspace{40 pt}if \segregation{$\kappa$}$\leq\lfloor\frac{n}{2}\rfloor$;

\hspace{60 pt}$D=$\constructdisc{$\kappa$};

\hspace{60 pt}$\gamma=$\constructcurve{$D$};

\hspace{60 pt}print $(\gamma, D)$;

\hspace{60 pt}STOP;

\hspace{40 pt}for each line of segregation $\ell$;

\hspace{60 pt}\GIRTH{$n$,$S-\ell$,$\dd$};

\end{ttfamily}

\subsection{Proof of Algorithm.}

\begin{claim}
The algorithm of Section~\ref{surface algorithm} works and is effective.
\end{claim}

\begin{proof}
The algorithm terminates after finite time because at each stage there 
are finitely many $\dd$-discs, finitely many $\dd$-half-discs, and 
finitely many lines of segregation.  There are finitely many stages in the
recursion because cutting along $\ell$ reduces $\ec{S}$ by one.

To show the algorithm works, it suffices to show that if $\girth(S)\leq n$,
the algorithm terminates with output.  This is vacuously true for $\chi(S)=2$.
So assume we have shown it for $\chi(S)\geq k+1$.  We will prove
it for $\chi(S)=k$.

Assume that $\girth(S)\leq n$.  
Let $\dd$ be a locally minimal disc system.
Let $\gamma$ be a closed curve in $S$ with girth $n$.  Let $D$ be a 
girth-realizing disc for $\gamma$.  Apply the Subordination Lemma to $D$.

There is either a $\dd$-disc subordinate to $D$, or two $\dd$-half-discs 
independently subordinate to $D$.  Depending on the case, let $\dint{a}$ 
be that $\dd$-disc, or the $\dd$-half-disc with smaller complexity.
Let $m$ be the complexity of $\dint{a}$.  Then $\gamma$ homotopically
intersects $\dint{a}$ exactly $n$ times.  If $\dint{a}$ is the $\dd$-disc
(resp. $\dd$-half-disc), then $m\leq n$ (resp. 
$m\leq\lfloor\frac{n}{2}\rfloor$).

By Lemma~\ref{segregated regions}, either \segregation{$\dint{a}$}$\leq m$, 
or $\gamma$ is isotopic to a curve disjoint from some line of segregation.  
Note that if $\chi (S)=0$, the first must be true.  In the former case, 
an application of \texttt{CONSTRUCT.DISC} and \texttt{CONSTRUCT.CURVE} 
finds a curve $\gamma$ and disc $D$ satisfying $|\gamma\cap\dd|\leq\max(1,n)$.
In the latter case, a girth-minimizing curve is contained 
in $S-\ell$ for some line of segregation 
(see the proof of Lemma~\ref{segregated regions}).
Then apply \texttt{GIRTH} to $S-\ell$, whose Euler characteristic 
is $\chi(S)=k+1$.
\end{proof}

\clearpage

\renewcommand{\baselinestretch}{1.6}\small\normalsize


\end{document}